%% file: main.tex
\documentclass{ACmod}

\input{defs}

\usepackage{tikz}
\usepackage[tocflat]{tocstyle}

\begin{document}

\input{title}

\input{abstract}

\maketitle

\setcounter{tocdepth}{1}
\tableofcontents

\input{intro}
\input{amodel}

\input{geomms}
\input{KZ}
\input{gauge}

\input{formalism}
\input{roadmap}
\input{vertices}
\input{computer}
\input{computation}

\input{biblio}
\end{document}

%% file: defs.tex
\DeclareFontFamily{U}{rsf}{}
\DeclareFontShape{U}{rsf}{m}{n}{
  <5> <6> rsfs5 <7> <8> <9> rsfs7 <10-> rsfs10}{}
\DeclareMathAlphabet{\mathscr}{U}{rsf}{m}{n}

\DeclareMathAlphabet{\mathgth}{U}{euf}{m}{n}

\DeclareFontFamily{U}{cyr}{}
\DeclareFontShape{U}{cyr}{m}{n}{
  <5> wncyr5 <6> wncyr6 <7> wncyr7 <8> wncyr8 <9> wncyr9 <10-> wncyr10}{}
\DeclareMathAlphabet{\mathcyr}{U}{cyr}{m}{n}

\input cyracc.def

\DeclareSymbolFont{bbold}{U}{bbold}{m}{n}
\DeclareSymbolFontAlphabet{\mathbbold}{bbold}

\makeatletter
\def\operator@font{\sf}
\makeatother

\newcommand{\hol}{\mathsf{hol}}

\newcommand{\cA}{{\mathscr A}}

\newcommand{\cW}{{\mathcal W}}

\newcommand{\cC}{{\mathscr C}}

\newcommand{\cF}{{\mathscr F}}

\newcommand{\cM}{{\mathscr M}}

\newcommand{\cO}{{\mathscr O}}

\newcommand{\tM}{\widetilde{M}}

\newcommand{\dzbar}{d\bar{z}}
\newcommand{\mi}{\mathtt{i}}
\newcommand{\twopii}{{2\pi\mi}}
\newcommand{\imtau}{\imag(\tau)}
\newcommand{\taubar}{\bar{\tau}}

\newcommand{\dg}{{\mathsf{dg}}}
\newcommand{\deform}{\mathsf{def}}
\newcommand{\abs}{{\mathsf{abs}}}
\newcommand{\framed}{\mathsf{fr}}
\newcommand{\cat}{{\mathsf{cat}}}
\newcommand{\geom}{{\mathsf{geom}}}

\newcommand{\PD}{{\mathsf{PD}}}

\newcommand{\D}{{\mathsf D}_{\mathsf{coh}}^b}

\newcommand{\chk}{{\scriptscriptstyle\vee}}

\newcommand{\dR}{\mathsf{dR}}

\newcommand{\pt}{{\mathsf{pt}}}

\newcommand{\res}{{\mathsf{res}}}
\newcommand{\hres}{{\mathsf{hres}}}

\newcommand{\bbone}{\mathbbold{1}}
\newcommand{\bbC}{\mathbb{C}}
\newcommand{\bbZ}{\mathbb{Z}}
\newcommand{\bbQ}{\mathbb{Q}}
\newcommand{\bbH}{\mathbb{H}}
\newcommand{\bbR}{\mathbb{R}}

\newcommand{\Kahler}{{\mathsf{K\ddot{a}hler}}}
\newcommand{\Mukai}{\mathsf{M}}
\newcommand{\Poincare}{{\mathsf P}}
\newcommand{\Tate}{\mathsf{Tate}}
\newcommand{\modular}{\mathsf{mod}}
\newcommand{\GM}{\mathsf{GM}}
\newcommand{\GGM}{\mathsf{GGM}}
\newcommand{\Dubrovin}{{\mathsf D}}
\newcommand{\complex}{{\mathsf{cx}}}

\newcommand{\sbt}{{\scalebox{0.5}{\textbullet}}}

\DeclareMathOperator{\Sym}{Sym}
\DeclareMathOperator{\Fuk}{Fuk}

\DeclareMathOperator{\Tr}{Tr}

\DeclareMathOperator{\KZ}{KZ}
\DeclareMathOperator{\wt}{wt}
\DeclareMathOperator{\HH}{HH}
\DeclareMathOperator{\HC}{HC}
\DeclareMathOperator{\HP}{HP}

\DeclareMathOperator{\SL}{SL}

\DeclareMathOperator{\imag}{Im}
\DeclareMathOperator{\id}{id}
\DeclareMathOperator{\ev}{ev}

\DeclareMathOperator{\Ext}{Ext}
\DeclareMathOperator{\RHom}{RHom}

\DeclareMathOperator{\KS}{KS}

\DeclareMathOperator{\Real}{Re}

\newcommand{\ra}{\rightarrow}

\newcommand{\lra}{\longrightarrow}
\newcommand{\contract}{\,\lrcorner\,}

\newcommand{\Z}{\mathbb{Z}}

\newcommand{\iso}{\cong}

\newcommand{\chkE}{\check{E\,}\!}
\newcommand{\chkX}{\check{X\,}\!}

\newcommand{\del}{\partial}
\newcommand{\delbar}{\bar{\partial}}

\renewcommand{\phi}{\varphi}

%% file: title.tex
\title{Computing a categorical Gromov-Witten invariant}

\author[C\u ald\u araru and Tu]{% 
Andrei C\u ald\u araru\footnote{Partially supported by the National
Science Foundation through grant number DMS-1200721.} and Junwu Tu}

\address{Andrei C\u ald\u araru, Mathematics Department,
University of Wisconsin--Madison, 480 Lincoln Drive, Madison, WI
53706--1388, USA.}
\address{Junwu Tu, Mathematics Department, University of Missouri,
Columbia, MO 65211, USA.}

%% file: abstract.tex
\begin{abstract} 
  {\sc Abstract:} We compute the $g=1$, $n=1$ B-model Gromov-Witten
  invariant of an elliptic curve $E$ directly from the derived
  category $\D(E)$.  More precisely, we carry out the computation of
  the categorical Gromov-Witten invariant defined by
  Costello using as target a cyclic $\cA_\infty$ model of
  $\D(E)$ described by Polishchuk.

  This is the first non-trivial computation of a positive genus
  categorical Gromov-Witten invariant, and the result agrees with the
  prediction of mirror symmetry: it matches the classical
  (non-categorical) Gromov-Witten invariants of a symplectic 2-torus
  computed by Dijkgraaf.
\end{abstract}

%% file: intro.tex
\section{Introduction}

\paragraph
The initial form of mirror symmetry, as described in 1991 by
Candelas-de la Ossa-Green-Parkes~\cite{COGP}, centered on the
surprising prediction that the genus zero Gromov-Witten invariants of
a quintic threefold $\chkX$ could be computed by solving a
differential equation governing the variation of Hodge structure
associated to another space, the so-called {\em mirror quintic} $X$.
Many other such mirror pairs $(X, \chkX)$ were later found in
physics, satisfying similar relationships between the genus zero
Gromov-Witten invariants of $\chkX$ and the variation of Hodge
structure of $X$.

\paragraph
A far-reaching generalization of mirror symmetry was proposed several
years later by Kontsevich~\cite{Kon} in his address to the 1994
International Congress of Mathematicians.  He conjectured that the
more fundamental relationship between the spaces $X$ and $\chkX$
in a mirror pair should be the existence of a derived equivalence
between the derived category $\D(X)$ of coherent sheaves on $X$ and
the Fukaya category $\Fuk(\chkX)$ of $\chkX$.  This statement
became known as the homological mirror symmetry conjecture.

\paragraph
Implicit in Kontsevich's proposal was the idea that the equality of
numerical invariants predicted by the original version of mirror
symmetry should follow tautologically from the homological mirror
symmetry conjecture.  To achieve this one needs to construct {\em
  categorical Gromov-Witten invariants}: invariants associated to an
(enhanced) triangulated category $\cC$, with the property that they
recover the classical Gromov-Witten invariants of the space
$\chkX$ when the target category $\cC$ is taken to be
$\Fuk(\chkX)$.  Once one has such invariants, evaluating them on
$\D(X)$ yields new invariants of $X$, the so-called {\em B-model
  Gromov-Witten invariants} of $X$.  These invariants are defined for
{\em any} genus, not just for genus zero.  (The genus zero B-model
invariants are expected to match the data of the variation of Hodge
structures used before.)  The categorical nature of the construction
automatically implies, for a pair of spaces $(X, \chkX)$ which
satisfies homological mirror symmetry, that the B-model invariants of
$X$ match the Gromov-Witten invariants of $\chkX$.

\paragraph
Genus zero categorical Gromov-Witten invariants satisfying the desired
properties were defined in 2015 by Ganatra-Perutz-Sheridan~\cite{GPS}
following ideas of Saito~\cite{Sai1, Sai2} and
Barannikov~\cite{Bar}.  However, according to the authors, this
approach does not extend to positive genus.

For arbitrary genus Costello~\cite{Cos} proposed a definition of
categorical invariants associated to a cyclic $\cA_\infty$ algebra (or
category), following ideas of Kontsevich and Soibelman~\cite{KonSoi}.
Unfortunately many details of~\cite{Cos} were left open, and computing
explicit examples turned out to be a difficult task.  Costello
(unpublished) computed one example where the target algebra is the
ground field (corresponding to the case where the target space $X$ is
a point).  No other explicit computations of Costello's invariants
exist. Costello-Li~\cite{CosLi} wrote in 2010:
\begin{quote} 
    {\em A candidate for the B-model partition function associated to a
    Calabi-Yau category was proposed
    in~\cite{Cos1},~\cite{Cos1},~\cite{KonSoi} based on a
    classification of a class of 2-dimensional topological field
    theories. Unfortunately, it is extremely difficult to compute this
    B-model partition function.}
\end{quote}

\paragraph
In this paper we compute the $g=1$, $n=1$ B-model categorical
invariant of an elliptic curve $E_\tau$, starting from Costello's
definition and using as input an $\cA_\infty$ model of the derived
category $\D(E_\tau)$ proposed by Polishchuk~\cite{Pol}.  It is the
first computation of a categorical Gromov-Witten invariant with
non-trivial target and positive genus.

More precisely, for a complex number $\tau$ in the upper-half plane
$\bbH$ let $E_\tau$ denote the elliptic curve of modular parameter
$\tau$, $E_\tau = \bbC/\bbZ\oplus\bbZ\tau$.
For each such $\tau$ we compute a complex number
$F_{1,1}^{\mathsf B}(\tau)$, the corresponding B-model categorical
invariant.  Regarding the result as a function of $\tau$ we obtain the
so-called {\em B-model Gromov-Witten potential}, a complex-valued
function on the upper half plane.

\paragraph
\label{subsec:mspred}
Mirror symmetry predicts the result of the above computation.  There
is a standard way to collect the classical $g=1$, $n=1$ Gromov-Witten
invariants of an elliptic curve in a generating power series, the {\em
  A-model Gromov-Witten potential} $F_{1,1}^{\mathsf A}(q)$.  The
result is known through work of Dijkgraaf~\cite{Dij}:
\[ F_{1,1}^{\mathsf A}(q) = -\frac{1}{24} E_2(q). \] 
Here $E_2$ is the standard Eisenstein holomorphic, quasi-modular 
form of weight 2, expanded at $q=\exp(2\pi\mi\tau)$.  For this
computation to give a non-trivial answer we need to insert at the one
puncture the class $[\pt]^\PD$ which is Poincar\'e dual to a point.

The prediction of mirror symmetry is that the A- and B-model
potentials should match after the K\"ahler and complex moduli spaces
are identified via the mirror map, which in the case of elliptic
curves takes the form
\[ q = \exp(2\pi\mi \tau). \]
Thus the prediction of mirror symmetry is that the B-model potential
should equal
\[ F_{1,1}^{\mathsf B}(\tau) = -\frac{1}{24} E_2(\tau). \]

\paragraph
\label{subsec:findxi}
To get our computation off the ground we need a cyclic
$\cA_\infty$-algebra model of the derived category $\D(E_\tau)$.  Such
an algebra was described by Polishchuk~\cite{Pol}, using structure
constants that are modular, almost holomorphic forms.  We will use both
Polishchuk's original algebra, and a gauge-equivalent modification of
it whose structure constants are quasi-modular, holomorphic forms.
The interplay between calculations in these two models, via the
Kaneko-Zagier theory of quasi-modular forms, will form a central part of
our final computation.

\paragraph
Like in the classical Gromov-Witten calculation, in order to get a
non-trivial answer in the B-model computation we need to insert a
certain Hochschild class $[\xi] \in HH_{-1}(A_\tau)$ at the puncture,
mirror dual to $[\pt]^\PD$.  This class will be represented by the
Hochschild chain in $A_\tau^{\otimes 1}$
\[ \xi = \frac{1}{\tau-\taubar} \dzbar. \]
(The identification of $E_\tau$ with $\bbC/\bbZ\oplus \bbZ\tau$ yields
a well-defined class $\dzbar$ in $H^1(E_\tau, \cO_{E_\tau})$.  This
group is a direct summand of $A_\tau$.  Therefore $\xi$ is a
well-defined element of homological degree $(-1)$ of the algebra
$A_\tau$, and as such it gives rise to a class in $HH_{-1}(A_\tau)$.)

The following theorem is the main result of this paper.

\begin{Theorem}
\label{thm:mainthm}
  With insertion the class $[\xi]$, Costello's categorical
  Gromov-Witten invariant of $A_\tau$ at $g=1$, $n=1$ equals
  \[ F_{1,1}^{\mathsf B}(\tau) = -\frac{1}{24} E_2(\tau). \] 
\end{Theorem}

\paragraph
We interpret this result in two ways.  On one hand we think of it as
confirmation of the mirror symmetry prediction at $g=1$ as
in~(\ref{subsec:mspred}).  On the other hand, through the prism of
homological mirror symmetry we can view our result as a statement about
the Fukaya category of the family $\chkE_\rho$ which is mirror to
the family $E_\tau$ of elliptic curves.  Indeed, by work of
Polishchuk-Zaslow~\cite{PolZas} we know that homological mirror
symmetry holds for elliptic curves.  The authors construct an equivalence
\[ \D(E_\tau) \iso \Fuk(\chkE_\rho), \] 
where $\chkE$ is the 2-torus mirror to $E_\tau$, endowed with a
certain complexified K\"ahler class $\rho$.  Therefore our computation,
which is {\em a priori} about $\D(E_\tau)$, can be reinterpreted as a
calculation about $\Fuk(\chkE_\rho)$.  From this perspective we
regard Theorem~\ref{thm:mainthm} as verification of the prediction
that Costello's categorical Gromov-Witten invariants of the Fukaya
category agree (in this case) with the classical ones of the
underlying space, as computed by Dijkgraaf.

\paragraph
There is one important aspect of Costello's work that we have
suppressed in the above discussion.  In order to extract an actual
Gromov-Witten potential from a cyclic $\cA_\infty$-algebra $A$ (as
opposed to a line in a certain Fock space) we need to choose a
splitting of the Hodge filtration on the periodic cyclic homology of
$A$.  The correct splitting is forced on us by mirror symmetry.  The
Hochschild and cyclic homology of $A_\tau$ agree with those of
$E_\tau$, as they are derived invariants.  Under this identification,
a splitting of the Hodge filtration is the choice of a splitting of
the natural projection
\[ H^1_\dR(E_\tau) \ra H^1(E_\tau, \cO_{E_\tau}) = HH_{-1}(A_\tau).  \] 
Mirror symmetry imposes the requirement that the lift of $[\xi]$ must
be invariant under monodromy around the cusp, which in turn uniquely
determines the lifting.  It is with this choice that we carry out the
computations in Theorem~\ref{thm:mainthm}.  See
Section~\ref{sec:computer} for more details. 

\paragraph
There is another approach to higher genus invariants in the B-model,
due to Costello and Li~\cite{CosLi, Li, Li3, Li2}, inspired by the
BCOV construction in physics~\cite{BCOV}.  

These other invariants also depend on a choice of splitting of the
Hodge filtration.  In their works Costello and Li analyzed the
BCOV-type invariants of elliptic curves obtained from {\em arbitrary}
splittings of the Hodge filtration and they showed that these
invariants satisfy the Virasoro constraints.  Moreover, they studied a
family of splittings depending on a parameter $\sigma\in \bbH$ and
they proved the modularity of the corresponding BCOV potentials.  The
monodromy invariant splitting that we consider corresponds to the
limiting splitting $\sigma \rightarrow \mi\infty$.  We have learned
the idea that this is the correct one for mirror symmetry from
conversations with Costello and Li.

\paragraph
The BCOV-type invariants have the advantage that they give a more
geometric definition of B-model Gromov-Witten invariants for
Calabi-Yau spaces, and are also more easily computed than the original
categorical ones of Costello~\cite{Cos}.  In fact, for elliptic curves
Li was able to establish mirror symmetry at arbitrary genus for BCOV
B-model invariants, and to directly compute the potential functions in
any genus and for arbitrary insertions.

However, the BCOV-type invariants are fundamentally different from the
ones we study in this paper, in that they are not {\em a priori}
categorical: knowing homological mirror symmetry does not allow one to
conclude the equality of the A- and B-model invariants.  Moreover, the
BCOV approach does not immediately generalize to other non-geometric
situations wherein one only has a category, and not an underlying space.

\paragraph {\bf Outline of the paper.}
Section~\ref{sec:A-model} outlines Dijkgraaf's computation in the
classical setting.  Section~\ref{sec:geomms} discusses mirror symmetry
in the geometric setting.  The next two sections review modular forms,
Kaneko-Zagier theory, and Polishchuk's $\cA_\infty$
algebra. Costello's general formalism is outlined in
Section~\ref{sec:formalism}, and the next section contains a roadmap
to the computation for elliptic curves.  Section~\ref{sec:vertices}
describes a computation, essentially due to Costello, of the string
vertices for $\chi=-1$.  The last two sections present two different
ways to compute the Gromov-Witten invariant $F_{1,1}^{\mathsf B}$ that
we want.  The first method involves reducing the problem to a very
large linear algebra computation, which is then solved by computer.
The second method, presented in Section~\ref{sec:computation}, gives a
purely mathematical deduction of the result, using a comparison
between computations in the holomorphic and modular gauges,
respectively.

\paragraph {\bf Standing assumptions.} We work over the field of
complex numbers $\bbC$.  Throughout the paper we will need to use
various comparison results between algebraic homology theories
(Hochschild, cyclic) and geometric ones (Hodge, de Rham).  Most of
these comparison results are in the literature; however, some appear
to be known only to the specialists but are not published.  In
particular we have not been able to find in the literature a
comparison between the algebraic Getzler-Gauss-Manin connection and
the classical geometric one.  We tacitly assume that they agree, but
this should be considered a conjectural result.  (A similar assumption
is made in~\cite{GPS}.)

\paragraph {\bf Acknowledgments.}  We would like to thank Nick
Sheridan, Kevin Costello, Si Li, Dima Arinkin, Alexander Polishchuk,
and Jie Zhou for patiently listening to the various problems we ran
into at different stages of the project, and for providing insight.
Stephen Wright's explanation of the use of $L^1$-optimization
techniques has been extremely helpful for computing a small solution
to the lifting problem.  Some of the larger computations were carried
out on the SBEL supercomputer at the University of Wisconsin, whose
graceful support we acknowledge.  

%% file: amodel.tex
\section{The classical invariants}
\label{sec:A-model}

In this section we outline Dijkgraaf's computation~\cite{Dij} of the
classical $g=1$, $n=1$ Gromov-Witten invariants of elliptic curves.

\paragraph
Let $\chkE = \bbR^2/\Z^2$ denote the two dimensional torus, endowed
with any complex structure making it into an elliptic curve.  (The
specific choice of complex structure will not matter.)  We upgrade
$\chkE$ to a symplectic manifold by choosing any symplectic form
$\omega$ in $H^2(\chkE, \bbR)$ whose area is one.

We also fix a point $P_0 \in \chkE$, which we think of as determining the
origin of the group structure of $\chkE$.

\paragraph
For any $\beta\in H_2(\chkE, \Z)$ the moduli space
$\overline{M}_{1,1}(\chkE, \beta)$ has virtual real dimension two.  If
$\beta = d\cdot[\chkE]$ for some integer $d>0$ then
\[ \overline{M}_{1,1}(\chkE, \beta) = M_{1,1}(\chkE, \beta) \] 
and the virtual dimension agrees with the actual dimension. In fact
the moduli space $M_{1,1}(\chkE,\beta)$ parametrizes in this case the
pairs $(f, P)$ where $f:E\ra \chkE$ is a degree $d$ isogeny onto
$\chkE$ from another elliptic curve $E$, and $P$ is any point on
$E$. (There are only finitely many such isogenies possible, therefore
$\dim_\bbR M_{1,1}(X, \beta) = 2$.)

\paragraph
\label{subsec:dijk}
To compute an actual numerical Gromov-Witten invariant we need to
insert a cohomology class at the marked point, so that the integrand
is a 2-form.  Inserting a $\psi$-class from $M_{1,1}$ gives zero, so
the only choice left is to pull-back a 2-form $\alpha$ from $\chkE = T^2$
via the evaluation map
\[ \ev_1: M_{1,1}(\chkE, \beta) \ra \chkE. \] 
The natural choice used is to take $\alpha = [\pt]^\PD$, the
Poincar\'e dual class to a point. If we think of this point as being
$P_0$, then the associated Gromov-Witten invariant for
$\beta=d\cdot[\chkE]$ will compute the number of isogenies
$f:E\ra \chkE$, up to isomorphism.  (Inserting $\alpha$ has the effect
of requiring $P$ to map to $P_0$.)  This number is well known -- it
equals the sum $(\sum_{k\,|\,d} k)$ of divisors of $d$.  This can be
seen by counting the number of matrices with integer coefficients, of
determinant $d$, up to $\SL(2, \Z)$ conjugation.

\paragraph
One encapsulates the result of this computation into the classical
Gromov-Witten potential of the A-model
\[ F_{1,1}^{\mathsf A}(q) = \sum_{\beta \in H_2(\chkE,\Z)} \langle
  \alpha \rangle_{1,1}^{\chkE, \beta} q^{\langle \beta, \omega\rangle}
  = -\frac{1}{24} + \sum_{d\geq 1} \sum_{k\,|\,d} k q^d. \] 
(The leading coefficient of $-1/24$ arises from a computation on
$\overline{M}_{1,1}$ which we will omit.)  This formula equals
$-\frac{1}{24}E_2(q)$, where $E_2$ is the Eisenstein holomorphic,
quasi-modular normalized form of weight 2; see Section~\ref{sec:KZ}
for details on modular forms.  Note that while the function $E_2$ is
not modular, it is still invariant under $\tau \mapsto \tau+1$.  Thus
the above formula is well defined, despite the integral ambiguity in
taking $\tau=\frac{1}{2\pi\mi}\log(q)$.

%% file: geomms.tex
\section{Geometric mirror symmetry}
\label{sec:geomms}

In this section we discuss the identifications that mirror symmetry
prescribes between structures of the A- and the B-models for elliptic
curves.  We place ourselves in the classical geometric context, where
one deals with spaces and not with categories.  

\paragraph
At its core mirror symmetry is an identification between two families
of geometric structures, the A-model and the B-model.  The A-model is
usually a trivial family of complex manifolds, endowed with a varying
{\em complexified K\"ahler class} (this notion is a generalization of
the usual K\"ahler class, see below).  The B-model family is a varying
family of complex manifolds.  The {\em mirror map} is an isomorphism
\[ \cM^\Kahler \iso \cM^\complex \] 
between the moduli space $\cM^\Kahler$ of complexified K\"ahler
classes and the moduli space $\cM^\complex$ of complex structures.
This isomorphism is not defined everywhere, but only in the
neighborhood of certain limit points of these spaces, the so-called
{\em large volume} and {\em large complex structure} limit points. 

\paragraph
In the case of elliptic curves both the A- and the B-model families
have descriptions in terms of familiar structures.  We begin by
describing the A-model family.

Let $\chkE$ denote the 2-torus $\bbR^2/\bbZ^2$
endowed with some complex structure, as in Section~\ref{sec:A-model}.
The moduli space of complexified K\"ahler structures on $\chkE$ is
defined to be
\[ \cM^\Kahler(\chkE) = \left (H^2(\chkE,\bbR)/H^2(\chkE,
    \bbZ) \right )
  \oplus \mi\cdot \left \{\omega\in H^2(\chkE, \bbR)~|~ \int_{\chkE}
  \omega > 0\right \}. \]
In other words a complexified K\"ahler class $\rho$ on $\chkE$ can
be written as $\rho = b+iA$ where $b$ is a form in
$H^2(\chkE,\bbR)/H^2(\chkE, \bbZ)$ and $A$ is a symplectic
form with positive area.

The moduli space $\cM^\Kahler$ is naturally isomorphic to $\bbH/\bbZ$.
The correspondence $\rho \leftrightarrow \tau$ thus identifies a
neighborhood of $\mi\cdot\infty$ in $\cM^\Kahler(\chkE)$ (the
large volume limit point) with a neighborhood of the cusp in the
moduli space of complex structures of elliptic curves (the large
complex structure limit point).  Indeed, only the subgroup $\bbZ$ of
translations by one of $\SL(2,\bbZ)$ acts in a small enough
neighborhood of the cusp.

The map $\rho \leftrightarrow \tau$ is the mirror map for
elliptic curves described by Polishchuk-Zaslow~\cite{PolZas}.  We have
used it in order to identify the A- and B-model potentials
in~(\ref{subsec:mspred}).

Due to the periodicity of $b$ and $\Real(\tau)$ and the positivity of $A$
and $\imtau$ it often makes more sense to use exponential coordinates
on $\cM^\Kahler$ and $\cM^\complex$. We will write
$q = \exp(2\pi\mi\rho)$ or $q=\exp(2\pi\mi\tau)$ depending on the
context, with the hope that this will not cause any confusion.

\paragraph
On the A-model side the trivial bundle over $\cM^\Kahler$ with fiber
$H^{1-\sbt}(\chkE, \bbC)$ carries the structure of a graded variation of
polarized semi-infinite Hodge structures (VSHS), in the sense
of~\cite{Bar},~\cite{GPS}.  We shall only need parts of this
structure: a graded fiberwise pairing (the Poincar\'e pairing
$\langle\,-\, ,\,-\,\rangle_\Poincare$ given by wedge and integrate)
and a flat connection, the {\em Dubrovin connection}
$\nabla^\Dubrovin$ (see~\cite{GPS}):
\[ \nabla_{q\del_q}(x) = q\del_q(x) - x\cup [\omega]. \] 

\paragraph
\label{subsec:Mukaipairing}
We have a similar structure of polarized VSHS over $\cM^\complex$
coming from the B-model.  The fiber over $\tau\in\cM^\complex$ of the
underlying vector bundle is the graded vector space
$HH_{\sbt-1}(E_\tau)$.  The fiberwise pairing is given by the {\em Mukai
  pairing} $\langle\,-\, ,\,-\,\rangle_\Mukai$.  The connection is the
Gauss-Manin connection $\nabla^\GM$ after identifying $HH_\sbt(E_\tau)$
with $H^\sbt_{\dR}(E_\tau)$ as $\bbZ/2\bbZ$-graded vector spaces. (We are
using the fact that the Hodge-de Rham spectral sequence degenerates on
$E_\tau$).

We remind the reader the formula for the Mukai pairing.  The
Hochschild-Kostant-Rosenberg isomorphism gives an identification
\[ HH_\sbt(E_\tau) \iso \bigoplus_{q-p=\sbt} H^p(E_\tau, \Omega^q_{E_\tau}).
\] 
The Mukai pairing then becomes
\[ \langle u, v \rangle_\Mukai = \frac{1}{\twopii}\int_{E_\tau}
  (-1)^{|u|} u \wedge v, \] 
where $|u| = q$ for $u \in H^p(E_\tau, \Omega^q_{E_\tau})$.  See
Ramadoss~\cite{Ram} for an explanation of this sign.  The factor of
$1/\twopii$ arises from the comparison of the Serre duality pairing
(the residue pairing) with the wedge-and-integrate pairing.  Note that
this factor is missing in the comparison Conjecture~1.14
in~\cite{GPS}.  Also note that the Todd class of an elliptic curve is
trivial, so there is no correction from it as in~\cite{CalHH2}.

\paragraph
The mirror map gives an isomorphism between the A-model polarized VSHS
and the B-model one.  The following theorem identifies the classes
$[\Omega]\in HH_1(E_\tau)$ and $[\xi]\in HH_{-1}(E_\tau)$ which
correspond to the classes $1\in H^0(\chkE_\rho)$ and
$[\pt]^\PD \in H^2(\chkE_\rho)$, respectively, under this isomorphism.

\begin{Proposition}
\label{prop:defOmegaxi}
Under the mirror map $\rho\leftrightarrow\tau$ the class
$1\in H^0(\chkE_\rho)$ corresponds to the class of the global
holomorphic volume form
\[ [\Omega] = [\twopii\cdot dz]\in HH_1(E_\tau). \] 
Similarly, the class $[\pt]^\PD \in H^2(\chkE_\rho)$ corresponds to
the class
\[ [\xi] = \frac{1}{\tau-\taubar} [\dzbar] \in
  HH_{-1}(E_\tau). \]
\end{Proposition}

\begin{Proof}
Because we are only interested in what happens in a neighborhood of
the large complex limit point we have well defined forms $dz$ and
$\dzbar$ which give bases of $H^0(E_\tau, \Omega^1_{E_\tau})$ and
$H^1(E_\tau \cO_{E_\tau})$, respectively, on each elliptic curve
$E_\tau$.  The classes $[\Omega]$ and $[\xi]$ are therefore pointwise
multiples of these forms.  The goal is to identify which multiples
they are.

In the A-model we have
\[ \langle 1, \nabla^\Dubrovin_{q\del_q}(1) \rangle_\Poincare = \langle 1,
  -[\omega] \rangle_\Poincare = -1. \]
The relationships $q = \exp(\twopii \rho)$ and
$\rho\leftrightarrow \tau$ force the identification
\[ q\del_q = \frac{1}{\twopii}\del_\tau. \]
Therefore we require in the B-model to have 
\[ \langle\Omega, \nabla^\GM_{\frac{1}{\twopii}\del_\tau} \Omega
  \rangle_\Mukai = -1. \] 
A straightforward calculation with periods shows that 
\[ \nabla_{\del_\tau}^\GM(dz) = \frac{1}{\tau-\taubar}(dz-\dzbar) \]
and this forces the class $[\Omega]$ to equal $[\twopii \cdot dz]$.

Similarly, the A-model identity
\[ \langle 1, [\pt]^\PD \rangle_\Poincare = 1 \]
forces the equality
\[ \langle [\Omega], [\xi] \rangle_\Mukai = 1,\] 
which in turn implies
\[ [\xi] = \frac{1}{\tau-\taubar} [\dzbar] \]
as stated in~(\ref{subsec:findxi}).
\end{Proof}

% \paragraph
% The choice of global holomorphic form $\Omega$ on $E_\tau$ also yields a
% trivialization of the canonical bundle of $E_\tau$. Thus the Serre
% duality pairing becomes a perfect pairing of degree one on the sheaf
% cohomology $H^*(E_\tau, \cF)$ of any coherent sheaf $\cF$ on $E_\tau$.
% For example, if we think of $\id_\cO$ and $\dzbar$ as basis vectors
% for $H^0(E_\tau, \cO_{E_\tau})$ and $H^1(E_\tau, \cO_{E_\tau})$,
% respectively, then this pairing gives
% \[ \langle \id_\cO, \dzbar \rangle = \frac{1}{\twopii} \int_{E_\tau} 1
%   \wedge \dzbar \wedge \Omega  = 2\mi \imtau. \]
% We will use this in Section~\cite{sec:ainfalgebra} to fix the correct
% cyclic pairing on the algebra $A_\tau$ (which will be different from
% the one used by Polishchuk in~\cite{Pol}).

\paragraph
We conclude this section with a discussion of the splitting of the
Hodge filtration on the de Rham cohomology $H^1_\dR(E_\tau)$.  For an
elliptic curve $E_\tau$ the Hodge filtration is expressed by the short
exact sequence
\[ 0 \ra H^0(E_\tau, \Omega^1_{E_\tau}) \ra H^1_\dR(E_\tau) \ra
  H^1(E_\tau, \cO_{E_\tau})\ra 0, \]
or, equivalently, by the sequence
\[ 0 \ra HH_1(E_\tau) \ra H^1_\dR(E_\tau) \ra HH_{-1}(E_\tau) \ra 0. \]
We have natural basis vectors $[dz]$ and $[\xi]$ of the first and last
vector spaces in this sequence.  Thus choosing a splitting of the Hodge
filtration means picking a lift $[\xi+f(\tau) dz]$ of $[\xi]$ from
\[ HH_{-1}(E_\tau) = H^1(E_\tau, \cO_{E_\tau}) \]
to $H^1_\dR(E_\tau)$ for every $\tau$.
\medskip

\noindent
Mirror symmetry prescribes that the correct lift must be {\em
  invariant} under monodromy around the cusp.  The following lemma
characterizes this lift.

\begin{Lemma}
\label{lem:geomsplitting}
The following conditions are equivalent for a family $[\tilde{\xi}]^\geom$
of lifts of the family of Hochschild classes $[\xi]$:
\begin{enumerate}
\item[(i)] The lift $[\tilde{\xi}]^\geom$ is invariant under monodromy around
  the cusp for all $\tau\in\bbH$.
\item[(ii)] The family $[\tilde{\xi}]^\geom$ is flat with respect to the
  Gauss-Manin connection.
\end{enumerate}
\end{Lemma}

\begin{Proof}
Fix an identification of $E_\tau$ with $\bbC/\bbZ\oplus\bbZ\tau$.
This determines a basis $A, B$ of cycles in $H_1(E_\tau)$
corresponding to the paths from 0 to 1 and from 0 to $\tau$, respectively.
Under the monodromy $\tau\mapsto \tau+1$ around the cusp the basis
$(A, B)$ maps to $(A, A+B)$.

Let $A^*, B^*$ denote the basis in $H^1_\dR(E_\tau)$ dual to the basis
$(A, B)$.  Then under the same monodromy the pair $(A^*, B^*)$ maps to
$(A^*-B^*, B^*)$ (the inverse transpose matrix).  We have
\begin{align*}
dz & = A^* + \tau B^* \\
\dzbar & = A^* + \taubar B^*.
\end{align*}
It follows that the invariant cocycle $B^*$ is expressed in the $(dz,
\dzbar)$ basis as 
\[ B^* = \frac{1}{\tau-\taubar}(dz-\dzbar). \]
For a class $[\tilde{\xi}] \in H^1_\dR(E_\tau)$ to be invariant under
monodromy around the cusp it must be a multiple of $B^*$.  And indeed
there exists a unique monodromy invariant lift of $[\xi]$ from
\[ HH_{-1}(E_\tau) = H^1_\dR(E_\tau)/HH_1(E_\tau) \]
to $H^1_\dR(E_\tau)$, namely
\[ -B^* = \frac{1}{\tau-\taubar}(\dzbar-dz). \]
It is now obvious that this family is $\nabla^\GM$-flat and this
condition uniquely identifies the family in (ii).
\end{Proof}

\paragraph
\label{subsec:KS}
We will also need in Section~\ref{sec:computation} the explicit form
of the Kodaira-Spencer class
\[ \KS(\del_\tau) = -\frac{1}{\tau-\taubar} \frac{\del}{\del z}\dzbar
  \in H^1(E_\tau, T_{E_\tau}). \]
This follows from the relationship of the Kodaira-Spencer class and the Gauss-Manin
connection, given by the formula
\[ \nabla^\GM_{\del_\tau}([\Omega]) \mod H^{1,0}(E_\tau) = \KS(\del_\tau) \contract
  [\Omega]. \]

%% file: KZ.tex
\section{Quasi-modular forms and Kaneko-Zagier theory}
\label{sec:KZ}

Before we can describe Polishchuk's $\cA_\infty$ model for the derived
category of an elliptic curve we need to review the theory of
quasi-modular, holomorphic forms and its relationship to the theory of
almost holomorphic, modular forms.  Their interplay is described by
Kaneko-Zagier theory, see~\cite[5.1]{KanZho} and~\cite{KanZag}.

\paragraph
The ring of holomorphic, modular forms for the group
$\Gamma=\SL(2,\bbZ)$ is isomorphic to $\bbC[E_4, E_6]$, where $E_4$
and $E_6$ are the Eisenstein modular forms of weights 4 and 6,
respectively.  (Since there are competing notations in the literature
we reserve the notation $E_k$ to mean the corresponding
normalized form, in the sense that the function has been rescaled to
satisfy $\lim_{\tau \ra \mi\cdot\infty} E_k(\tau) = 1$.)

\paragraph
For the purposes of this paper we need to consider certain functions
on $\bbH$ which are not at the same time holomorphic and modular.  Of
particular interest to us is the Eisenstein form $E_2$: it is still
holomorphic, but it does not satisfy the usual transformation law with
respect to $\Gamma$, see~\cite[4.1]{KanZho}.  We will call elements of
the ring $\tM(\Gamma) = \bbC[E_2, E_4, E_6]$ {\em quasi-modular}
holomorphic forms; the weight of such a form is defined by declaring
the weight of $E_{2k}$ to be $2k$.

\paragraph
\label{subsec:defe2star}
While the form $E_2$ is not modular, the following simple modification
\[ E_2^*(\tau) =  E_2(\tau)-\frac{3}{\pi^2}\frac{\twopii}{\tau-\taubar} \]
is -- it satisfies the same transformation rule as regular modular
forms of weight two with respect to the group $\Gamma$.  However,
$E_2^*$ is no longer holomorphic.  The functions in the ring
$\cO(\bbH)[\frac{1}{\tau-\taubar}]$ which are modular in the
above sense (for arbitrary weights) and satisfy suitable growth
condition at the cusp form a ring $\widehat{M}(\Gamma)$.
Its elements will be called {\em almost holomorphic} modular forms.
One can show that
\[ \widehat{M}(\Gamma) = \bbC[E_2^*, E_4, E_6]. \]

\paragraph
\label{subsec:delhat}
The rings $\tM(\Gamma)$ and $\widehat{M}(\Gamma)$ are closed under the
actions of certain differential operators.   The former is closed
under the action of $\del_\tau$, while the latter is closed under
\[ \widehat{\del}_\tau = \del_\tau + \frac{\mathsf{wt}}{\tau-\taubar}, \]
where $\mathsf{wt}$ denotes the weight of the form on which
$\widehat{\del}_\tau$ acts.  

The theorem below is the main result of Kaneko-Zagier theory~\cite{KanZag}.

\begin{Theorem}
\label{thm:KZ}
The operator 
\[ \phi : \widehat{M}(\Gamma) \ra \tM(\Gamma) \]
which maps an almost holomorphic, modular form
\[ F(\tau, \taubar) = \sum_{m\geq 0} \frac{F_m(\tau)}{(\tau-\taubar)^m} \]
to its ``constant term'' part $F_0(\tau)$ is a differential ring
isomorphism.  Its inverse will be denoted by 
\[ \KZ: \tM(\Gamma) \ra \widehat{M}(\Gamma). \]
\end{Theorem}

%% file: gauge.tex
\section{Polishchuk's algebra and its holomorphic modification}
\label{sec:gauge}

In this section, we review Polishchuk's cyclic $\cA_\infty$ algebra
$A_\tau$ associated to an elliptic curve $E_\tau$.  This presentation
uses almost holomorphic, modular forms.  The main result of this
section asserts that the $\cA_\infty$ algebra $A_\tau^\hol$ obtained
by replacing these structure constants by their images under the
homomorphism $\phi$ of Theorem~\ref{thm:KZ} yields an algebra which is
{\em gauge equivalent} to $A_\tau$ and varies holomorphically with
$\tau$.

\paragraph
For $\tau\in \bbH$ let $E_\tau$ be the corresponding elliptic
curve.  We will denote by $P_0$ its origin.  The derived category
$\D(E_\tau)$ of coherent sheaves on the elliptic curves $E_\tau$ is
compactly generated by the object $\cF = \cO\oplus L$, where $L$ is
the degree one line bundle $L = \cO(P_0)$.  Therefore the dg-algebra
\[ A_\tau^\dg = \RHom_{E_\tau}^\sbt(\cF, \cF) \] 
is derived equivalent to $E_\tau$. By homological perturbation we can
transfer the dg-algebra structure from $A_\tau^\dg$ to a
quasi-equivalent $\cA_\infty$ algebra structure on
\[ A_\tau = H^\sbt(A_\tau^\dg) = \Ext^\sbt_{E_\tau}(\cF, \cF). \] 
The graded vector space $A_\tau$ is concentrated in cohomological
degrees zero and one, and 
\[ \dim A_\tau^0 = \dim A_\tau^1 = 3. \]

\paragraph
Polishchuk~\cite{Pol} explicitly computes an $\cA_\infty$ structure on
$A_\tau$ using a particular choice of homotopy between $A_\tau$ and
$A_\tau^\dg$.  In order to express his result Polishchuk chooses a
basis for the 6-dimensional algebra $A_\tau$ consisting of:
\begin{itemize}
\item The identity morphisms $\id_\cO:\cO\ra \cO$ and $\id_L:L\ra L$
  in $A_\tau^0$. 

\item The morphisms $\xi :\cO \ra \cO[1]$ and $\xi_L: L \ra L[1]$ in
  $A_\tau^1$ given by $\xi = \dzbar$ and similarly for $\xi_L$.

\item The degree zero morphism $\theta:\cO\ra L$ in $A_\tau^0$ given
  by the theta function $\theta = \theta(z, \tau)$.

\item An explicit dual morphism $\eta:L \ra \cO[1]$ in $A_\tau^1$
  given by the formula in~\cite[2.2]{Pol}.
\end{itemize}

\paragraph
\label{subsec:basis}
For our purposes it will be more convenient to work with a slightly
modified basis for $A_\tau$, obtained by dividing the elements of
degree one in Polishchuk's basis ($\eta$, $\xi$ and $\xi_L$) by the
factor $\tau-\taubar$.  We leave the degree zero elements unchanged.
The rescaled elements will still be denoted by $\eta$, $\xi$, $\xi_L$;
they differ by a factor of $1/\twopii$ from the ones that Polishchuk
considers in~\cite[2.5, Remark 2]{Pol}).

With this rescaling the new element $\xi$ agrees with the expression
in~(\ref{subsec:findxi}).

\paragraph
\label{subsec:ainf}
At this point we need to decide which sign conventions to use for the
$\cA_\infty$ axioms.  There are (at least) two such conventions in the
literature.  In this paper we will use the so-called {\em shifted}
signs, which are more common in the symplectic geometry literature and
are easier to use.  See~\cite{Cho} for a discussion of these two
conventions and of ways to translate between them.  Note, however,
that Polishchuk is using in~\cite{Pol} the unshifted sign convention;
we will make the change to the shifted convention as we go.
Following~\cite{Cho} we will denote by $|x|'$ the shifted degree of an
element $x\in A_\tau$, $|x|' = |x|-1$.

\paragraph
Once we have fixed a global holomorphic volume form $\Omega$ on
$E_\tau$ as in Proposition~\ref{prop:defOmegaxi} we get an induced
symmetric non-degenerate pairing on $A_\tau$ arising from Serre
duality.  Correcting it by the shifted sign conventions we obtain
a skew-symmetric pairing on $A_\tau$ given by
\begin{align*}
\langle \xi, \id_\cO \rangle = -1,  & & \langle \xi_L, \id_L \rangle = -1,  \\
\langle \id_\cO, \xi \rangle = 1,  & & \langle \id_L, \xi_L \rangle = 1,  \\
\langle \theta, \eta \rangle = 1, & & \langle \eta, \theta \rangle = -1.
\end{align*}
This pairing will give the {\em cyclic} structure on $A_\tau$. 

\paragraph
We will now describe Polishchuk's formulas for the multiplications
\[ \mu_k: A[1]^{\otimes k} \ra A[2]. \]
These operations are non-trivial only when $k$ is even.  The first
product $\mu_2$ is the usual Yoneda product on the $\Ext$ algebra
$A_\tau$, with a sign adjustment: the only non-zero products (beside
the obvious ones involving $\id$ and $\id_L$) are $\theta\eta = -\xi$
and $\eta\theta = \xi_L$.  (We will use reverse notation for
composition of morphisms, so that $\theta\eta$ means
$\eta\circ \theta$.  The same convention shall be used for the higher
products, keeping in line with Polishchuk's conventions.)

\paragraph
\label{subsec:cycainf}
An important feature of Polishchuk's formulas is that they respect
cyclic symmetry with respect to the inner product on $A_\tau$.  We
will say that a homogeneous map (of some degree)
\[ c_k:A^{\otimes k} \ra \bbC \]
is cyclic if it is invariant under cyclic rotation of its arguments,
up to a Koszul sign determined by the shifted degrees of its
arguments.   The products $\mu_k$ will be cyclic in the sense that the
tensors
\[ c_k(x_1, \ldots, x_k) = \langle \mu_{k-1}(x_1, \ldots, x_{k-1}),
  x_k \rangle \]
are cyclic in the above sense.  An $\cA_\infty$ algebra endowed with a
pairing of degree $d$ which makes  the maps $c_k$ cyclic will be said
to be Calabi-Yau of degree $d$, or CY$[d]$, or cyclic of degree $d$.

\paragraph
All the higher multiplications of $A_\tau$ can be deduced from a
single one using this cyclic symmetry.  The complete list is presented
in~\cite[Theorem 2.5.1]{Pol}, but for reference we list the formula
for only one such multiplication.  Explicitly, let $a$, $b$, $c$, $d$
be non-negative integers, and let $s = a+b+c+d$.  Then for $s$ odd we
have
\[ \mu_k(\xi_L^a,\eta,\xi^b,\theta,\xi_L^c,\eta,\xi^d) = 
(-1)^{a+b+{ {s+1}\choose{2}}} \frac{1}{a!b!c!d!} \cdot
\frac{1}{(\twopii)^{s+1}}\cdot g_{a+c,b+d}\cdot \eta, \] 
and zero otherwise. In this formula $k =s+3$, and  
\[ g_{a+c,b+d}\in \widehat{M}(\Gamma)_{s+1} \] 
is a certain almost holomorphic modular form of weight $s+1$
defined by Polishchuk in~\cite[1.1]{Pol}.

\paragraph
\label{subsec:wtbasis}
We associate a notion of weight to the basis elements of $A_\tau$, as
follows:
\begin{eqnarray*}
\wt(\id) & = \wt(\id_L) & = 0, \\
\wt(\theta) & = \wt(\eta) & = 1/2, \\
\wt(\xi) & = \wt(\xi_L) & = 1.
\end{eqnarray*}
Note that with respect to these assignments of weights the
multiplications $\mu_k$ are of total weight zero.

\paragraph
Even though the dependence on $\tau\in \bbH$ of the family of elliptic
curves $E_\tau$ is holomorphic, the associated family of $\cA_\infty$
algebras $A_\tau$ does not depend holomorphically on $\tau$: the
structure constants $g_{a+c, b+d}$ are not holomorphic.  We think of
the multiplications $\{\mu_k\}$ as giving a 2-cocycle $\mu^*$ in the
Hochschild cochain $CC^\sbt(A_\tau)$.  (The differential in this
complex is $[\mu^*,\,-\,]$.)  The failure of holomorphicity is measured
by the 2-cochain $\delbar\mu^*$.

The following result shows that the anti-holomorphic dependence of
this family is nevertheless trivial to first order (i.e., it is zero
in Hochschild cohomology). The proof is a straight-forward
computation.

\begin{Theorem}
\label{thm:exactness}
Let $CC^\sbt(A_\tau)$ denote the Hochschild cochain complex of $A_\tau$,
with differential $[\mu^*, \,-\,]$.  Let $\psi_3: A_\tau[1]^{\otimes
  3} \ra A_\tau[1]$ be the cyclic Hochschild 1-cochain defined by
\[ \psi_3(\xi,\theta,\xi_L)=-\frac{1}{\twopii}\cdot
  \frac{1}{(\tau-\taubar)^2} \cdot \theta \] 
and its cyclic rotations with respect to the pairing.  Then we have
\[ \delbar \mu^* = [\mu^*, \psi_3]. \]
\end{Theorem}

\paragraph
\label{subsec:defatauhol}
The deformation induced by $\delbar\mu$ is trivial not
just to first order: the bounding element $\psi_3$ can be integrated
to an $\cA_\infty$ quasi-isomorphism
\[ f : A_\tau \ra A_\tau^\hol,\] 
from $A_\tau$ to a {\em holomorphic} family $A_\tau^\hol$ of
$\cA_\infty$ algebras over $\mathbb{H}$. The process is similar to
integrating pseudo-isotopies of $\cA_\infty$ algebras,
see~\cite[Proposition 9.1]{Fuk}.  We will now present the construction
of the algebras $A_\tau^\hol$ and of the $\cA_\infty$ quasi-isomorphism
$f$.

The holomorphic family $A_\tau^{\hol}$ has the same underlying vector
space as $A_\tau$ and the same basis.  Its structure constants are
obtained simply by applying the map $\phi$ of Theorem~\ref{thm:KZ} to
the structure constants of $A_\tau$, so that for example, the
multiplications of $A_\tau^\hol$ are obtained by cyclic symmetry from
\[ \mu_k^\hol(\xi_L^a,\eta,\xi^b,\theta,\xi_L^c,\eta,\xi^d) = 
(-1)^{a+b+{ {s+1}\choose{2}}} \frac{1}{a!b!c!d!} \cdot
\frac{1}{(\twopii)^{s+1}}\cdot \phi( g_{a+c,b+d})\cdot \eta. \] 
Since $\phi$ is a ring map, the maps $\{\mu^\hol_k\}$ also satisfy the
$\cA_\infty$ relations. 

\paragraph
One can explicitly compute $\phi(g_{a,b})$ using the recurrence
formulas in~\cite[Proposition 2.6.1 (ii)]{Pol}. For example we have
\begin{align*}
\phi(g_{1,0})&=\phi(e_2^*) = \phi(e_2-\frac{\twopii}{\tau-\taubar})=e_2,\\
\phi(g_{2,1})&=\phi(-g_{1,0}^2+\frac{5}{6}g_{3,0})=-e_2^2+\frac{5}{6}\cdot
               3!\cdot e_4=-e_2^2+5e_4. 
\end{align*}
In general $g_{a,b}$ is a polynomial expression with rational
coefficients in $e_{2k}^*$; $\phi(g_{a,b})$ is obtained by replacing
$e_2^*$ by $e_2$ in the same expression, leaving all the other
terms unchanged.

The following theorem describes recursively the construction of the
maps $f_n$ which assemble to give a quasi-equivalence $A_\tau \iso
A_\tau^\hol$.  The proof is a straightforward (though tedious)
inductive calculation. 

\begin{Theorem}
\label{thm:holgauge}
Inductively define multi-linear maps
\[ f_n: A_\tau[1]^{\otimes n} \ra A^\hol_\tau[1] \]
by setting $f_1=\id$ and
\[ f_n:=\sum_{i\geq 1, j\geq 1, k\geq 1, i+j+k=n} \int
  \psi_3(f_i\otimes f_j \otimes f_k) \, d\taubar. \] 
Here the integration symbol $\int (\ldots)\,\, d\taubar$ is
formally applied to the coefficients of the tensors, and is defined by
\[ \int \frac{1}{(\tau-\taubar)^m} d\taubar:=
  \frac{1}{(m-1)}\frac{1}{(\tau-\taubar)^{m-1}} \mbox{ for } m\geq 2.\]  
Then the maps $\{f_n\}_{n\geq 1}$ form a cyclic $\cA_\infty$ quasi-isomorphism
\[ f: A_\tau\ra A^\hol_\tau.\]
\end{Theorem}

%% file: formalism.tex
\section{Costello's formalism}
\label{sec:formalism}

In this section we review Costello's general definition of
categorical Gromov-Witten invariants. In the next section we will
specialize this construction to the case of the $g=1$, $n=1$ invariant
of an elliptic curve.

\paragraph
Let $A$ be a cyclic $\cA_\infty$ algebra whose pairing is of degree
$d$ as in~(\ref{subsec:cycainf}), and let $\HH = \HH_\sbt(A)[d]$ denote
the shifted Hochschild homology of $A$.  The categorical Gromov-Witten
potential that we ultimately want to construct is an element
\[ F^\cat \in \Sym(u^{-1} \HH[u^{-1}])\series{\lambda}. \]  
Here $\lambda$ is a formal variable used to keep track of the Euler
characteristics of the moduli spaces and $u$ is the usual formal
variable of degree $-2$ which appears in the definition of
cyclic homology.  It is a placeholder for keeping track of
insertions of $\psi$ classes.

The invariant $F^\cat$ depends not only on the algebra $A$, but also
on a further choice: a splitting of the Hodge filtration on the
periodic cyclic homology $\HP_\sbt(A)$ of $A$.  The construction of
$F^\cat$ proceeds in two stages.  First one constructs an abstract
invariant $F^\abs$ which only depends on $A$. It is a state in the
homology of a certain Fock space associated to $A$. The choice of
splitting of the Hodge filtration is then used to identify the
homology of the Fock space with
$\Sym(u^{-1} \HH[u^{-1}])\series{\lambda}$ and hence to extract
$F^\cat$ from $F^\abs$.

\paragraph
\label{subsec:invariance}
The potential thus constructed is invariant under cyclic quasi-equivalence.
If $A'$ is another cyclic $\cA_\infty$ algebra and $f$ is a cyclic
quasi-equivalence $A\iso A'$, then $f$ induces an isomorphism
$\HP_\sbt(A) \iso \HP_\sbt(A')$ of filtered vector spaces.  The
original splitting $s$ of the Hodge filtration on $\HP_\sbt(A)$ then
determines a splitting $s'$ of the Hodge filtration on
$\HP_\sbt(A')$.  The invariance of the Gromov-Witten invariants means
that the potential associated to $(A,s)$ is the same as the one
associated to $(A',s')$ under the obvious identifications.

\paragraph
Individual Gromov-Witten invariants can be read off from $F^\cat$ as
coefficients of its power series expansion. For example, in the next
section we will want to compute the $g=1$, $n=1$ invariant of
Polishchuk's algebra $A_\tau$. The Euler characteristic of the curves
in $M_{1,1}$ is $\chi = -1$, thus we are interested in the
$\lambda^{-\chi} = \lambda^1$ part of $F^\cat$. We want to insert the
class $[\xi]$ (and no $\psi$-classes). This means that the
Gromov-Witten invariant we are interested in is the coefficient of
$[\Omega]u^{-1} \lambda$ in the power series expansion of
$F^\cat$: here $[\Omega]$ is the class in $HH_1(A_\tau)$ which is
Mukai dual to $[\xi]$ (Proposition~\ref{prop:defOmegaxi}), and we have
used $u^{-1}$ to denote no $\psi$-class insertions (one $\psi$ class
would give $u^{-2}$, etc.).

We begin by describing the structures on the cyclic chain complex of
$A$ that will be needed for the construction of $F^\abs$.

\paragraph {\bf The circle action on Hochschild chains.}
Let $V$ denote the shifted Hochschild cochain complex of $A$,
\[ V = (CC_\sbt(A), b)[d], \] 
whose homology is $\HH$.  Connes' $B$ operator is a degree one
operator on $V$.  It gives a homological circle action on the dg
vector space $V$.  We use it to form the periodic cyclic complex
\begin{align*}
V_\Tate & = \left (CC_\sbt(A)\laurent{u}, b+uB \right )[d] 
\intertext{and its negative cyclic subcomplex}
V^{hS^1} & = \left (CC_\sbt(A)\series{u}, b+uB \right )[d].
\end{align*} 
The homology of these complexes are the (shifted by d) periodic and
negative cyclic homology of $A$, respectively.

\paragraph {\bf The Mukai pairing.}
The Mukai pairing~\cite{CalMukai} is a non-degenerate symmetric
pairing of degree $2d$ on $\HH$.  Shklyarov~\cite{Shk} and
Sheridan~\cite{She} construct a lift of the Mukai pairing to the chain
level.  It is a symmetric bilinear form of degree $2d$ on $V$,
\[ \langle\,-\,,\,-\,\rangle_\Mukai : V\otimes V \ra \bbC[2d]. \] 
We will discuss this pairing further in~(\ref{subsec:action}).

\paragraph {\bf The higher residue pairing.}
The chain-level Mukai pairing on $V$ induces a $\bbC\laurent{u}$
sesquilinear pairing of degree $2d$ on $V_\Tate$,
\[ \langle\,-\, ,\,-\,\rangle_\hres: V_\Tate \times V_\Tate \ra
  \bbC\laurent{u}[2d], \] 
the {\em higher residue pairing}.  (See~\cite[Definition 2 and Section
5.42]{She} for details.)  For $\alpha, \beta\in CC_\sbt(A)$, $i, 
j\in \bbZ$ the higher residue pairing is given by 
\[ \langle \alpha u^i, \beta u^j \rangle_\hres = (-1)^i \langle
  \alpha, \beta \rangle_\Mukai \cdot u^{i+j}. \] 
Laurent power series only have finitely many negative powers of $u$,
hence the above formula extends by sesquilinearity to all of
$V_\Tate$. 

The higher residue pairing is a chain map with respect to $b+uB$,
which makes it descend to a pairing on periodic cyclic homology. The
space $V_\Tate$ is filtered by powers of $u$, and the pairing respects
this filtration. In particular, for $x, y\in \HC_\sbt^-(A)$ we have
$\langle x, y \rangle_\hres \in \bbC\series{u}$.

The constant term of the power series $\langle x, y \rangle_\hres$
equals the Mukai pairing
$\langle \overline{x}, \overline{y} \rangle_\Mukai$ of the reductions
of $x$ and $y$ to Hochschild homology.  However, the power series
$\langle x, y \rangle_\hres$ can have higher powers of $u$ which are
not computed by the above homology Mukai pairing.

\paragraph {\bf The ordinary residue pairing.}
We will also consider the {\em residue pairing}, the $\bbC$-valued
pairing of degree $2d-2$ on $V_\Tate$ obtained by taking the
coefficient of $u^{-1}$ (the residue) in the higher residue pairing.
It is skew-symmetric, and restricts to zero on the dg subspace
$V^{hS^1}$.  We will think of $V^{hS^1}$ as a Lagrangian subspace of
the symplectic vector space $V_\Tate$.

\paragraph {\bf Weyl algebra and Fock module.}
We associate a Weyl algebra $\cW$ to the symplectic vector space
$V_\Tate$ and a Fock space $\cF$ to its Lagrangian subspace
$V^{hS^1}$.  Explicitly, $\cW$ is defined as
\[ \cW = T^\sbt(V_\Tate)/([x,y] - \langle x, y \rangle_\res) \] 
where $[x,y]$ is the graded commutator of $x, y\in T^\sbt(V_\Tate)$, and
$\cF$ is the quotient of $\cW$ by the left-ideal generated by
$V^{hS^1}\subset V^\Tate$. It is an irreducible left $\cW$-module.

\paragraph {\bf String vertices.}
Sen-Zwiebach~\cite{SenZwi} note that there are certain distinguished
singular chains $S_{g,n}$, the {\em string vertices}, on the
uncompactified moduli spaces of curves $M_{g,n}$. They play the role
of approximations of the fundamental classes of the compactified
spaces $\overline{M}_{g,n}$.  The string vertices can be defined by
recursive relations encapsulated in a differential equation known as
the {\em quantum master equation}, see~\cite[Theorem 1]{Cos}.  We will
discuss string vertices in Section~\ref{sec:vertices}.

Moreover, Kontsevich and Soibelman argue that chains on moduli spaces
of curves act on the cyclic chain complex of any Calabi-Yau
$\cA_\infty$ algebra.  This action will be further explained 
in~(\ref{subsec:action}). 

\paragraph {\bf Deformed Fock module.}
\label{subsec:defrel}
Costello's main observation in~\cite{Cos} is that the action of the
string vertices $S_{g,n}$ on cyclic chains can be used to deform the
module structure on $\cF$.  The result is a
$\cW\series{\lambda}$-module $\cF^\deform$ which is a deformation of
$\cF$ over $\bbC\series{\lambda}$.

The deformed module $\cF^\deform$ is constructed as follows.  The
standard Fock module $\cF$ is the quotient of $\cW$ by the left ideal
generated by the relations
\[ \alpha u^k = 0\mbox{ for all } \alpha\in V, k\geq 0. \]
The module $\cF^\deform$ is obtained by taking the quotient of $\cW$
by the left ideal $I^\deform$ generated by the deformed relations
\[ \alpha u^k = \sum_{g,n} \rho(\iota(S_{g,n}))(\alpha u^k)
  \lambda^{2g-2+n} \]
for all $\alpha\in V$, $k\geq 0$.  Here $\iota$ is a certain operator
on chains that will be described in Section~\ref{sec:vertices}, and
$\rho$ denotes the action of singular chains on cyclic chains.

Costello argues that for reasonable $\cA_\infty$ algebras, for which
the Hodge filtration splits in the sense of~(\ref{def:splitHodge})
below, the homology $H_\sbt(\cF^\deform)$ is a flat deformation of the
Fock module $H_\sbt(\cF)$ over the classical (non-dg) Weyl algebra
$H_\sbt(\cW)$.

\paragraph
Fock modules over classical Weyl algebras are rigid, and thus the
above deformation must be trivial at the level of homology.  This
implies that there is an isomorphism
\[ H_\sbt(\cF^\deform) \stackrel{\sim}{\lra} H_\sbt(\cF)\series{\lambda} \]
of $H_\sbt(\cW)\series{\lambda}$-modules, unique up to multiplication by a
power series in $\bbC\series{\lambda}$ which begins with $1$.  

By its construction the module $H_\sbt(\cF^\deform)$ has a canonical
generator $\bbone$ (it is a quotient of $H_\sbt(\cW\series{\lambda})$
by a certain ideal).  Costello defines the abstract Gromov-Witten
potential
\[ F^\abs \in H_\sbt(\cF)\series{\lambda} \]
to be the image of $\bbone$ under the above isomorphism.  The
ambiguity from the above power series will not play a role for us
because we only care about the leading term of $F^\abs$ (the
$\lambda^1$ term).

\paragraph {\bf Splitting of the Hodge filtration.}
In order to extract a concrete categorical Gromov-Witten potential
\[ F^\cat \in \Sym(u^{-1} \HH[u^{-1}])\series{\lambda} \] 
from the abstract one we need to choose a splitting of the Hodge
filtration, a notion we now make precise.

Endow the graded vector space $\HH$ with a trivial circle operator and
define $\HH_\Tate = \HH\laurent{u}$.  The homology level Mukai pairing
on $\HH$ induces a higher residue pairing on $\HH_\Tate$ defined by
the same formula as the higher residue pairing.  The higher residue
pairing on $\HH_\Tate$ respects the {\em grading} with respect to
powers of $u$.  This grading induces a decreasing filtration on $\HH_\Tate$,
similar to the one on $H_\sbt(V_\Tate)$.

\begin{Definition}
\label{def:splitHodge}
A splitting of the Hodge filtration on the periodic cyclic homology of
$A$ is an isomorphism of filtered vector spaces 
\[ \HH_\Tate \stackrel{\sim}{\lra} H_\sbt(V_\Tate) \]
which respects the higher residue pairings. 
\end{Definition}
\medskip

\paragraph
Choosing a splitting of the Hodge filtration is equivalent to
assigning to each $x\in \HH_i(A)$ a lift $\tilde{x} \in \HC^-_i(A)$.
This assignment is required to satisfy the property that for $x, y\in
\HH_\sbt(A)$ we have 
\[ \langle \tilde{x}, \tilde{y} \rangle_\hres = \langle x,
  y\rangle_\Mukai. \] 
In other words we impose the condition that the higher residue pairing
evaluated on lifts of elements in $\HH$ must have no higher powers of
$u$.

\paragraph
Endowing the graded vector space $\HH_\Tate$ with the residue
Mukai pairing makes it into a symplectic vector space, and
$\HH^{hS^1}$ is a Lagrangian subspace.  This yields a corresponding Weyl
algebra $\cW_H$ and Fock space $\cF_H$. The latter can be identified
with $\Sym (u^{-1}\HH[u^{-1}])$ as a $\cW_H$-module via the action on
the generator $\bbone$: negative powers of $u$ act by multiplication,
non-negative powers of $u$ act by differentiation.

A choice of splitting of the Hodge filtration induces an isomorphism
of Weyl algebras
\[ \cW_H \stackrel{\sim}{\lra} H_\sbt(\cW) \]
and a corresponding isomorphism of Fock modules
\[ S:\Sym(u^{-1}\HH[u^{-1}])\series{\lambda} \iso \cF_H
  \stackrel{\sim}{\lra} H_\sbt(\cF)\series{\lambda}.\] 
The categorical Gromov-Witten potential $F^\cat$ that we ultimately
want is defined to be 
\[ F^\cat = S^{-1}(F^\abs), \]
the preimage under $S$ of the abstract Gromov-Witten potential.

%% file: roadmap.tex
\section{A roadmap to the computation}
\label{sec:roadmap}

We specialize Costello's formalism described in the previous section
to the computation of the $g=1$, $n=1$ invariant of elliptic curves.
We complete the main calculation of Theorem~\ref{thm:mainthm},
assuming a few intermediate computations which we put off for
Sections~\ref{sec:vertices}-\ref{sec:computation}.

\paragraph
Fix a complex number $\tau\in \bbH$ and consider the cyclic
$\cA_\infty$ algebra $A_\tau$ discussed in Section~\ref{sec:gauge},
whose pairing has degree $d=1$. We are interested in computing the
coefficient of $[\Omega] u^{-1} \lambda$ in the potential $F^\cat$
obtained from $A_\tau$, using a splitting of the Hodge filtration on
the periodic cyclic homology of $A_\tau$ which matches at the
geometric level the one described in Lemma~\ref{lem:geomsplitting}.

As explained above the computation proceeds in three steps:
\begin{itemize}
\item[(i)] We first compute the string vertices $S_{0,3}$ and
  $S_{1,1}$. (These are the only pairs $(g,n)$ with Euler
  characteristic $-1$.)  The answer will be expressed as certain
  linear combinations of ribbon graphs, using a model for chains on
  decorated moduli spaces of curves due to
  Kontsevich-Soibelman~\cite{KonSoi}.
\item[(ii)] We then find an explicit splitting of the Hodge filtration
  on $\HP_\sbt(A)$ by choosing a particular lift
  $[\tilde{\xi}]\in\HC_{-1}^-(A)$ of the class $[\xi]\in\HH_{-1}(A)$.
  This splitting is chosen to match the geometric one in
  Lemma~\ref{lem:geomsplitting}.
\item[(iii)] We complete the computation by combining the results of
  (i) and (ii) into a calculation with Fock modules.
\end{itemize}
In this section we will outline these steps, leaving the details of
the explicit computations for 
Sections~\ref{sec:vertices}-\ref{sec:computation}. 

\paragraph {\bf The string vertices.}
\label{subsec:strvert}
We want to represent singular chains on $M_{g,n}$ by linear
combinations of ribbon graphs, using classical results by
Strebel~\cite{Str} and Penner~\cite{Pen} that compare the ribbon graph
complex and the singular chain complex $C_\sbt(M_{g,n})$. For
technical reasons explained below the classical ribbon graph complex
is not flexible enough and we need a generalization described by
Kontsevich and Soibelman in~\cite[Section 11.6]{KonSoi}.  The details
will be discussed in Section~\ref{sec:vertices}.

The main reason this generalization is needed is so that we can write
a combinatorial model for the quantum master equation, which at the
level of singular chains was described in~\cite{Cos}.  The
combinatorial equation can then be solved recursively, degree by
degree in $\lambda$.  The first few terms of the solution are as
follows:
\begin{align*}
\iota(S_{0,3}) & =\frac{1}{2}\;\;\;
  \begin{tikzpicture}[baseline={([yshift=-2ex]current bounding
      box.center)},scale=0.3] 
\draw [thick] (0,0) to (0,2);
\draw [thick] (-0.2, 1.8) to (0.2, 2.2);
\draw [thick] (0.2, 1.8) to (-0.2, 2.2);
\draw [thick] (0,0) to (-2,0);
\draw [thick] (0,0) to (2,0);
\draw [thick] (-2.2,0) circle [radius=0.2];
\draw [thick] (2.2,0) circle [radius=0.2];
\node at (1.3, 2.3) {$u^{-1}$};
\node at (-2.2, 1) {$u^{-1}$};
\node at (2.2, 1) {$u^{-1}$};
\end{tikzpicture}
\intertext{and}
\iota(S_{1,1}) & =\frac{1}{24}\;\;\;
  \begin{tikzpicture}[baseline={(current bounding
      box.center)},scale=0.3] 
\draw [thick] (0,2) circle [radius=2];
\draw [thick] (-2,2) to (-0.6,2);
\draw [thick] (-0.8,2.2) to (-0.4,1.8);
\draw [thick] (-0.8, 1.8) to (-0.4, 2.2);
\draw [thick] (-1.4142, 0.5858) to [out=40, in=140] (1, 0.5);
\draw [thick] (1.25, 0.3) to [out=-45, in=225] (2, 0.2);
\draw [thick] (2,0.2) to [out=45, in=-50] (1.732, 1);
\node at (0.7, 2.2) {$u^{-2}$};
\end{tikzpicture}
\;\;-\;\;\frac{1}{4}\;\;\; 
\begin{tikzpicture}[baseline={(current bounding box.center)},scale=0.3]
\draw [thick] (0,2) circle [radius=2];
\draw [thick] (0,0) to (0,1.4);
\draw [thick] (-0.2, 1.2) to (0.2, 1.6);
\draw [thick] (-0.2, 1.6) to (0.2, 1.2);
\draw [thick] (0,0) to [out=80, in=180] (0.5, 1);
\draw [thick] (0.5,1) to [out=0, in=100] (0.9, 0.4);
\draw [thick] (0,0) to [out=-80, in=180] (0.5, -1);
\draw [thick] (0.5, -1) to [out=0, in=-100] (0.9, 0);
\node at (0,2.2) {$u^{-1}$};
\end{tikzpicture}.
\end{align*}
(The operator $\iota$ which turns an output into an input will also be
discussed in Section~\ref{sec:vertices}.)

\paragraph {\bf The correct lift.}
The class $[\Omega] \in \HH_1(A_\tau)$ has a canonical lift to cyclic
homology, therefore the only information needed to specify a splitting
of the Hodge filtration for $A_\tau$ is to choose a lift $[\tilde{\xi}]
\in \HC_{-1}^-(A_\tau)$ of the Hochschild homology class $[\xi]\in
\HH_{-1}(A_\tau)$ defined in Proposition~\ref{prop:defOmegaxi}.  Such
a lift will have the form 
\[ \tilde{\xi} = \xi + \alpha \cdot u + O(u^2) \]
for some chain $\alpha\in CC_1(A_\tau)$ that satisfies
\[ b(\alpha) = -B(\xi) = -1 \otimes \xi.\] 
For degree reasons any choice of $\alpha$ satisfying this condition
will uniquely determine a lift $\tilde{\xi}$: the higher order terms
can always be filled in to satisfy $(b+uB)(\tilde{\xi}) = 0$, and
these further choices do not change the class $[\tilde{\xi}]$ in
$\HP_{-1}(A_\tau)$.

We will argue in Section~\ref{sec:computer} that the correct lift
is characterized by the system of equations
\[
 \left \{ \begin{array}{l}
b(\alpha)  = -1\otimes \xi \\
b^{1|1}(\del_\tau\mu^*~|~ \alpha) = 0.
\end{array} \right .
\]
Here $b^{1|1}$ is the operator defined in~\cite[3.11]{She} and
$\del_\tau \mu^*$ is the derivative of the multiplications $\{\mu_k\}$
of $A_\tau$, viewed as a Hochschild 2-cocycle in $\HH^2(A_\tau)$.

These equations are linear in $\alpha$; we will argue in
Section~\ref{sec:computer} that large-scale computer calculations
allow us to solve this system explicitly. Alternatively, in the same
section we will argue by theoretical means that this system admits
a solution, and this solution is unique up to a $b$-exact term.

\paragraph {\bf The action of singular chains on cyclic chains.} 
\label{subsec:action}
Kontsevich-Soibelman construct an action of their generalized ribbon
graphs (thought of as chains on moduli spaces of curves) on the Hochschild
chain complex of an $\cA_\infty$ algebra $A$.  Specifically, to a
ribbon graph $\Gamma$ with $m$ inputs and $n$ outputs they 
associate a map 
\[ \rho(\Gamma) : CC_\sbt(A)^{\otimes m} \ra CC_\sbt(A)^{\otimes n}. \] 
of a certain degree depending on $\Gamma$. For example the chain-level
Mukai pairing is the degree $2d$ map associated to the following graph
of genus zero, two inputs, no outputs:
\[\begin{tikzpicture}[baseline={(current bounding box.center)},scale=0.3]
\draw [thick] (0,2) circle [radius=2];
\draw [thick] (-2,2) to (-0.6,2);
\draw [thick] (-0.8,2.2) to (-0.4,1.8);
\draw [thick] (-0.8, 1.8) to (-0.4, 2.2);
\draw [thick] (2,2) to (3.4,2);
\draw [thick] (3.2,2.2) to (3.6,1.8);
\draw [thick] (3.2, 1.8) to (3.6, 2.2);
\end{tikzpicture}. \]
Similarly, the coproduct on the shifted Hochschild homology of a
cyclic $\cA_\infty$ algebra (induced by the product on Hochschild
cohomology, under the Calabi-Yau identification of cohomology with
shifted homology) is given by the graph with one input and two outputs
below:
\[   \begin{tikzpicture}[baseline={([yshift=-2ex]current bounding
      box.center)},scale=0.3] 
\draw [thick] (0,0) to (0,2);
\draw [thick] (-0.2, 1.8) to (0.2, 2.2);
\draw [thick] (0.2, 1.8) to (-0.2, 2.2);
\draw [thick] (0,0) to (-2,0);
\draw [thick] (0,0) to (2,0);
\draw [thick] (-2.2,0) circle [radius=0.2];
\draw [thick] (2.2,0) circle [radius=0.2];
\end{tikzpicture}. \]

\paragraph {\bf The Fock space computation.}
We are interested in the coefficient of $[\Omega] u^{-1}\lambda^1$ in
the series
\[ F^\cat \in \Sym(u^{-1}\HH[u^{-1}])\series{\lambda}. \] 
The way the algebra $\cW_H$ acts this coefficient can be recovered as
the constant coefficient of $\lambda^1$ in the expansion of
$[\xi] \cdot F^\cat$, where we have denoted by $\cdot$ the action of
$\cW_H$. 

Since $F^\cat$ is the preimage of the generator $\bbone$ of
$H_\sbt(\cF^\deform)$ under the isomorphism of $\cW_H$-modules
\[ S: \Sym(u^{-1}\HH[u^{-1}])\series{\lambda} \iso
  H_\sbt(\cF)\series{\lambda}  \stackrel{\sim}{\lra}
  H_\sbt(\cF^\deform), \]
the element $[\xi]\cdot F^\cat$ is the unique element in
$\Sym(u^{-1}\HH[u^{-1}])\series{\lambda}$ satisfying
\[ S([\xi]\cdot F^\cat) = [\tilde{\xi}]\cdot \bbone =
  [\tilde{\xi}]. \]
(We have used here the fact that the isomorphism $\cW_H \iso
H_\sbt(\cW)$ induced by the splitting maps $[\xi]\in \cW_H$ to
$[\tilde{\xi}]$, and the second $\cdot$ above refers to the action of
$H_\sbt(\cW)$.) 

\paragraph
Write the lift $[\tilde{\xi}]$ of $[\xi]$ as
\[ \tilde{\xi} = \xi + \alpha u + \beta u^2 \] 
by collecting in $\beta u^2$ the terms with powers of $u$ greater than
2.  Consider $\tilde{\xi} = \tilde{\xi}u^0$ as a chain level element
in $\cF^\deform$.  In this module the deformed
relations~(\ref{subsec:defrel}) give the following equality
\[ \tilde{\xi}  = \tilde{\xi}u^0 = \left ( \tilde{T}_1(\tilde{\xi}) +
    \tilde{T}_2(\tilde{\xi}) + \tilde{T}_3(\tilde{\xi}) \right )\lambda +
  O(\lambda^2), \]
where $T_1, T_2, T_3$ are the results of inserting $\tilde{\xi}$ into
the three ribbon graphs that appear in the expressions for the string vertices $S_{0,3}$
and $S_{1,1}$ as in~(\ref{subsec:strvert}):
\begin{eqnarray*}
\tilde{T}_1(\tilde{\xi})  & = \ \ \ \dfrac{1}{2}\!\!\!\! & \rho\big(
              \begin{tikzpicture}[baseline={([yshift=-2ex]current
                  bounding box.center)},scale=0.3] 
\draw [thick] (0,0) to (0,2);
\draw [thick] (-0.2, 1.8) to (0.2, 2.2);
\draw [thick] (0.2, 1.8) to (-0.2, 2.2);
\draw [thick] (0,0) to (-2,0);
\draw [thick] (0,0) to (2,0);
\draw [thick] (-2.2,0) circle [radius=0.2];
\draw [thick] (2.2,0) circle [radius=0.2];
\node at (1.3, 2.2) {$u^{-1}$};
\node at (-2.2, 1) {$u^{-1}$};
\node at (2.2, 1) {$u^{-1}$};
\end{tikzpicture}\big) (\tilde{\xi})  = \ \ \,
\dfrac{1}{2}  \ \ \rho\big(
              \begin{tikzpicture}[baseline={([yshift=-2ex]current
                  bounding box.center)},scale=0.3] 
\draw [thick] (0,0) to (0,2);
\draw [thick] (-0.2, 1.8) to (0.2, 2.2);
\draw [thick] (0.2, 1.8) to (-0.2, 2.2);
\draw [thick] (0,0) to (-2,0);
\draw [thick] (0,0) to (2,0);
\draw [thick] (-2.2,0) circle [radius=0.2];
\draw [thick] (2.2,0) circle [radius=0.2];
\node at (-2.2, 1) {$u^{-1}$};
\node at (2.2, 1) {$u^{-1}$};
\end{tikzpicture}\big) (\xi) \\\ \\
\tilde{T}_2(\tilde{\xi}) & = \ \ \dfrac{1}{24} \!\!\!\! & \rho \big(\begin{tikzpicture}[baseline={(current
                  bounding box.center)},scale=0.3] 
\draw [thick] (0,2) circle [radius=2];
\draw [thick] (-2,2) to (-0.6,2);
\draw [thick] (-0.8,2.2) to (-0.4,1.8);
\draw [thick] (-0.8, 1.8) to (-0.4, 2.2);
\draw [thick] (-1.4142, 0.5858) to [out=40, in=140] (1, 0.5);
\draw [thick] (1.25, 0.3) to [out=-45, in=225] (2, 0.2);
\draw [thick] (2,0.2) to [out=45, in=-50] (1.732, 1);
\node at (0.7, 2.2) {$u^{-2}$};
\end{tikzpicture}\big) (\tilde{\xi}) \quad\quad \ =\ \,
\dfrac{1}{24} \ \rho \big(\begin{tikzpicture}[baseline={(current
                  bounding box.center)},scale=0.3] 
\draw [thick] (0,2) circle [radius=2];
\draw [thick] (-2,2) to (-0.6,2);
\draw [thick] (-0.8,2.2) to (-0.4,1.8);
\draw [thick] (-0.8, 1.8) to (-0.4, 2.2);
\draw [thick] (-1.4142, 0.5858) to [out=40, in=140] (1, 0.5);
\draw [thick] (1.25, 0.3) to [out=-45, in=225] (2, 0.2);
\draw [thick] (2,0.2) to [out=45, in=-50] (1.732, 1);
\end{tikzpicture}\big) (\alpha)\\\ \\
\tilde{T}_3(\tilde{\xi}) & =  \ \ - \dfrac{1}{4}\!\!\!\! & \rho\big( 
\begin{tikzpicture}[baseline={(current bounding box.center)},scale=0.3]
\draw [thick] (0,2) circle [radius=2];
\draw [thick] (0,0) to (0,1.4);
\draw [thick] (-0.2, 1.2) to (0.2, 1.6);
\draw [thick] (-0.2, 1.6) to (0.2, 1.2);
\draw [thick] (0,0) to [out=80, in=180] (0.5, 1);
\draw [thick] (0.5,1) to [out=0, in=100] (0.9, 0.4);
\draw [thick] (0,0) to [out=-80, in=180] (0.5, -1);
\draw [thick] (0.5, -1) to [out=0, in=-100] (0.9, 0);
\node at (0,2) {$u^{-1}$};
\end{tikzpicture}\big) (\tilde{\xi}) \quad\quad\ \, = 
- \dfrac{1}{4}\ \rho\big( 
\begin{tikzpicture}[baseline={(current bounding box.center)},scale=0.3]
\draw [thick] (0,2) circle [radius=2];
\draw [thick] (0,0) to (0,1.4);
\draw [thick] (-0.2, 1.2) to (0.2, 1.6);
\draw [thick] (-0.2, 1.6) to (0.2, 1.2);
\draw [thick] (0,0) to [out=80, in=180] (0.5, 1);
\draw [thick] (0.5,1) to [out=0, in=100] (0.9, 0.4);
\draw [thick] (0,0) to [out=-80, in=180] (0.5, -1);
\draw [thick] (0.5, -1) to [out=0, in=-100] (0.9, 0);
\end{tikzpicture}\big) (\xi).
\end{eqnarray*}
Since $\tilde{T}_2(\tilde{\xi})$ only depends on the chain $\alpha$, we will
sometimes denote it by $T_2(\alpha)$. 

\paragraph
Having explicit formulas for the operations $\mu_k$ of the algebra
$A_\tau$ and for the chain $\alpha$ allows us to compute the values of
the three expressions above.  This will be done with a computer
calculation in Section~\ref{sec:computer}, and the result is that
$T_2(\tilde{\xi}) = T_3(\tilde{\xi}) = 0$ and
\[ T_1(\tilde{\xi}) = \frac{1}{2} (\xi u^{-1})(\xi u^{-1}) \in \cW. \]
Thus in $\cF^\deform$ we have the equality
\[ \tilde{\xi} = \frac{1}{2} (\xi u^{-1})(\xi u^{-1}) \lambda +
    O(\lambda^2). \]

\paragraph
As written above the expression on the right hand side is not
explicitly in the image of $S$, because the classes $\xi$ are not
$b+uB$-closed while the classes in the image of $S$ always are.
However, we know abstractly that this must be the case, because the
map $S$ is an isomorphism of dg modules.  So all we need to do is
to rewrite the expression $ (\xi u^{-1})(\xi u^{-1})$ using the Weyl
algebra relations to make it appear as the image under $S$ of an
element in $\cW_H$.

Note that $(\tilde{\xi} u^{-1})(\tilde{\xi} u^{-1})$ is the image
under $S$ of $([\xi] u^{-1})([\xi] u^{-1})$.  For degree reasons the
terms $\alpha$ and $\beta u$ are in the left ideal
$I^\deform \mod \lambda^2$ (inserting them into any of the ribbon
graphs in $S_{0,3}$ and $S_{1,1}$ gives zero). Thus we have the
equality in $\cF^\deform$
\begin{align*}
(\xi u^{-1})(\xi u^{-1}) & = (\tilde{\xi}u^{-1} - \alpha - \beta u)
                           (\tilde{\xi}u^{-1} - \alpha - \beta u) \\ 
& = (\tilde{\xi}u^{-1})(\tilde{\xi}u^{-1}) - (\alpha+\beta
  u)(\tilde{\xi}u^{-1})\\ 
& = (\tilde{\xi}u^{-1})(\tilde{\xi}u^{-1}) - \langle \alpha, \xi
  \rangle_\Mukai.
\end{align*}
Here we have used the Weyl algebra relation 
\[ (\alpha+\beta u)(\tilde{\xi}u^{-1}) = (\tilde{\xi}u^{-1})
  (\alpha+\beta u) + \langle \alpha+\beta u, \tilde{\xi}u^{-1}
  \rangle_\res = \langle \alpha, \xi \rangle_\Mukai. \]

\paragraph
We put together all the previous computations to get 
\[ [\xi]\cdot F^\cat = S^{-1}(\tilde{\xi}) = \frac{1}{2} ([\xi] u^{-1})([\xi]
    u^{-1}) \lambda -\frac{1}{2} \langle \alpha, \xi \rangle_\Mukai
    \lambda + O(\lambda^2), \]
and therefore the categorical Gromov-Witten invariant we want (the
constant coefficient of $\lambda$ in $[\xi]\cdot F^\cat$) is  
\[ F_{1,1}^{\mathsf B}(\tau) = -\frac{1}{2} \langle \alpha, \xi
  \rangle_\Mukai. \] 
We will argue in Section~\ref{sec:computer} that
\[ \langle \alpha, \xi \rangle_\Mukai = \frac{1}{12} E_2(\tau), \] 
thus completing the computation of Theorem~\ref{thm:mainthm}:
\[ F_{1,1}^{\mathsf B}(\tau) = -\frac{1}{24} E_2(\tau). \]

%% file: vertices.tex
\section{String vertices}
\label{sec:vertices}
In this section we discuss the quantum master equation on chains on
moduli spaces of curves, and the construction of string vertices as
linear combinations of ribbon graphs with rational coefficients.

\paragraph
Let $M_{g,0,n}^\framed$ denote the coarse moduli space of smooth,
genus $g$ curves with $n$ framed (parametrized), ordered
outgoing boundaries.  The wreath product group
$(S^1)^n \ltimes \Sigma_n$ acts on the space $M_{g, 0,n}^\framed$,
where the circles $(S^1)^n$ act by rotation of the framing, and the
symmetric group $\Sigma_n$ acts by permutation of the
boundaries. Denote by
$M_{g,0, n}=M_{g,0, n}^\framed/((S^1)^n\ltimes \Sigma_n)$ the quotient
space, parametrizing smooth, genus $g$ curves with $n$
unframed, unordered outgoing boundaries. Let $C_\sbt(M_{g,0,n})$ be the
singular chain complex of $M_{g,0,n}$ with rational coefficients.

There is a degree one operator 
\[ \Delta: C_\sbt(M_{g,0,n}) \rightarrow C_{\sbt+1}(M_{g+1, 0, n-2}),\]
defined as the sum of all the possible ways of sewing pairs of
boundary components, with a full $S^1$ twist.  Twist-sewing also
defines an odd bracket
\[ \left\{ -, -\right\} : C_\sbt(M_{g_1,n_1})\otimes C_\sbt(M_{g_2,n_2})
  \rightarrow C_{\sbt+1}(M_{g_1+g_2,n_1+n_2-2}),\] 
where the gluing is performed by choosing one boundary from the first
surface and one from the second surface.

\begin{Theorem}
\label{thm:vertex}
(Costello~\cite[Section 4, Theorem 1]{Cos}) There exists, for each
pair of integers $g\geq 0$, $n>0$ with $2-2g-n<0$, a chain
$S_{g,n} \in C_\sbt(M_{g,n})$ of degree $6g-6+2n$ such that
\begin{itemize}
\item[(i)] $S_{0,3}$ is a $0$-chain of degree $\frac{1}{3!}$.
\item[(ii)] For all $g, n \geq 0$ and $2-2g-n>0$ we have
\begin{equation*}
\partial S_{g,n} + \sum_{\substack{g_1+g_2=g\\ n_1+n_2=n+2}}
\frac{1}{2}\left\{ S_{g_1,n_1}, S_{g_2,n_2}\right\} + \Delta
S_{g-1,n+2}=0.
\end{equation*}
\end{itemize}
Furthermore, such a collection $\left\{ S_{g,n}\right\}$ is unique up
to homotopy in the sense of~\cite[Section 8]{Cos}.
\end{Theorem}
\medskip

\paragraph
The equation in part (ii) of the theorem is a form of the {\em quantum
  master equation} for chains on moduli spaces, see~\cite{Cos}.  It
was first discovered by Sen and Zwiebach~\cite{SenZwi}.  The singular
chains $S_{g,n}$ solving the quantum master equation are called {\em
  string vertices}.

\paragraph
In order to perform explicit computations we need to use a
combinatorial model for chains on moduli spaces of curves
for which the twist-sewing operations make sense.  Such a
combinatorial model is described by Kontsevich and
Soibelman~\cite[Section 11.6]{KonSoi}.  Their construction gives a
chain model $C_\sbt^{{\sf comb}}(M_{g,m,n}^\framed)$ for the moduli space
$M_{g,m,n}^\framed$ of genus $g$ curves with $m$ inputs, $n$
outputs which are ordered and framed.  The moduli space
$M_{g,m,n} = M_{g,m,n}^\framed/\left ((S^1)^n\ltimes S_{m,n}\right )$ is defined as
before.

However, the ribbon graph model only exists for $m\geq 1$ -- ribbon
graphs must have at least one incoming boundary component. To remedy
this Costello notes that there is a map 
\[ \iota : C_\sbt( M_{g,0,n}) \rightarrow C_\sbt(M_{g,1,n-1}),\]
obtained by switching the designation of one of the outgoing
boundaries in $M_{g,0, n}$ to ``incoming'' and summing over all such
choices of outgoing boundary.  

\paragraph
We can describe the chains $S_{g,n}$ combinatorially, via the map
$\iota$.  Taking coinvariants of the $(S^1)^n\ltimes \Sigma_{1,n-1}$ action
on $C_\sbt^{{\sf comb}}(M_{g,1,n-1}^\framed)$ produces a combinatorial model
of $C_\sbt(M_{g,1,n-1})$ of the form
\[ C_\sbt^{{\sf comb}}(M_{g,1,n-1})= \Big( (u_1^{-1},\ldots,u_n^{-1})
  C_\sbt^{{\sf comb}}(M_{g,1,n-1}^\framed)[u_1^{-1},\dots,u_n^{-1}]
  \Big)_{\Sigma_{1,n-1}}.\] 
More details can be found in~\cite[Section 5]{Cos}. 

Thus, a chain in $C_\sbt^{{\sf comb}}(M_{g,1,n-1})$ is a decorated ribbon
graph whose input and outputs are labeled by negative powers of
$u$. For example, the chains $\iota(S_{0,3})$ and $\iota(S_{1,1})$
that we have described in~(\ref{subsec:strvert}) are elements in
$C_\sbt(M_{0,1,2})$ and $C_\sbt(M_{1,1,0})$, respectively.  The vertices
labeled with a cross are input vertices, which for
Kontsevich-Soibelman are in the set $V_{\mathsf{in}}$. The white
vertices (little circles) in $\iota(S_{0,3})$ denote output vertices
in $V_{\mathsf {out}}$.

\paragraph
We now explain where the coefficients that appear in the chains
defined in~(\ref{subsec:strvert}) come from.  For 
\begin{align*}
\iota(S_{0,3}) & =\frac{1}{2}\;\;\;
  \begin{tikzpicture}[baseline={([yshift=-2ex]current bounding
      box.center)},scale=0.3] 
\draw [thick] (0,0) to (0,2);
\draw [thick] (-0.2, 1.8) to (0.2, 2.2);
\draw [thick] (0.2, 1.8) to (-0.2, 2.2);
\draw [thick] (0,0) to (-2,0);
\draw [thick] (0,0) to (2,0);
\draw [thick] (-2.2,0) circle [radius=0.2];
\draw [thick] (2.2,0) circle [radius=0.2];
\node at (1.3, 2.3) {$u^{-1}$};
\node at (-2.2, 1) {$u^{-1}$};
\node at (2.2, 1) {$u^{-1}$};
\end{tikzpicture}
\end{align*}
the coefficient is $\frac{1}{2}$ instead of $\frac{1}{3!}$ because
when applying the map $\iota$ there are three choices of outgoing
boundary to turn to an incoming one.

Next we compute
$\iota\big(\Delta (S_{0,3}) \big)\in C^{{\sf comb}}_\sbt(M_{1,1,0})$
in the moduli space of genus one curves with one unparametrized input
boundary component, and zero outgoing boundary components. Since
${ 3 \choose 2}=3$, the chain $\Delta(S_{0,3})$ is
$\frac{1}{2} = 3\cdot\frac{1}{6}$ times the result of self-sewing of any
two boundary components. Self-sewing of $2$ outgoing boundary
components can be described combinatorially:
\begin{equation*} \iota\big( \Delta (S_{0,3})\big) = \frac{1}{2} \,\,\begin{tikzpicture}[baseline={(current bounding box.center)},scale=0.3]
\draw [thick] (0,2) circle [radius=2];
\draw [thick] (-2,2) to (-0.6,2);
\draw [thick] (-0.8,2.2) to (-0.4,1.8);
\draw [thick] (-0.8, 1.8) to (-0.4, 2.2);
\draw [thick] (0,0) to [out=80, in=180] (0.5, 1);
\draw [thick] (0.5,1) to [out=0, in=100] (0.9, 0.4);
\draw [thick] (0,0) to [out=-80, in=180] (0.5, -1);
\draw [thick] (0.5, -1) to [out=0, in=-100] (0.9, 0);
\node at (0.7,2.2) {$u^{-1}$};
\end{tikzpicture}\,\, \in C_1^{{\sf comb}}(M_{1,1,0}).\end{equation*}

\paragraph
The string vertex $\iota(S_{1,1})$, being a degree $2$ chain in
$C^{{\sf comb}}_2(M_{1,0,1})$, must be a linear combination of the
form
\[ \iota(S_{1,1})=a \;\;\cdot\;\; \begin{tikzpicture}[baseline={(current bounding box.center)},scale=0.3]
\draw [thick] (0,2) circle [radius=2];
\draw [thick] (-2,2) to (-0.6,2);
\draw [thick] (-0.8,2.2) to (-0.4,1.8);
\draw [thick] (-0.8, 1.8) to (-0.4, 2.2);
\draw [thick] (-1.4142, 0.5858) to [out=40, in=140] (1, 0.5);
\draw [thick] (1.25, 0.3) to [out=-45, in=225] (2, 0.2);
\draw [thick] (2,0.2) to [out=45, in=-50] (1.732, 1);
\node at (0.7, 2.3) {$u^{-2}$};
\end{tikzpicture}
\;\;+\;\; b\;\;\cdot\;\;
\begin{tikzpicture}[baseline={(current bounding box.center)},scale=0.3]
\draw [thick] (0,2) circle [radius=2];
\draw [thick] (0,0) to (0,1.4);
\draw [thick] (-0.2, 1.2) to (0.2, 1.6);
\draw [thick] (-0.2, 1.6) to (0.2, 1.2);
\draw [thick] (0,0) to [out=80, in=180] (0.5, 1);
\draw [thick] (0.5,1) to [out=0, in=100] (0.9, 0.4);
\draw [thick] (0,0) to [out=-80, in=180] (0.5, -1);
\draw [thick] (0.5, -1) to [out=0, in=-100] (0.9, 0);
\node at (0,2.2) {$u^{-1}$};
\end{tikzpicture}\] 
for some $a,b\in \bbQ$.  (These are all the possible ribbon graphs
of even degree in $C_\sbt^{\mathsf{comb}}(M_{1,1,0})$.)  

Denote by $\partial$ the boundary operator on
$C^{{\sf comb}}_\sbt(M_{1,0,1})$, and by $D$ the circle operator which is
of degree one. The total differential on the equivariant
complex $C^{{\sf comb}}_\sbt(M_{1,1,0})$ is given by $(\partial+uD)$.
Then the quantum master equation for $n=1$, $g=1$ becomes
\[(\partial + u D) S_{1,1} + \Delta(S_{0,3}) = 0.\]
One can check that 
\begin{align*}
 \partial \Bigg ( \begin{tikzpicture}[baseline={(current bounding box.center)},scale=0.3]
\draw [thick] (0,2) circle [radius=2];
\draw [thick] (0,0) to (0,1.4);
\draw [thick] (-0.2, 1.2) to (0.2, 1.6);
\draw [thick] (-0.2, 1.6) to (0.2, 1.2);
\draw [thick] (0,0) to [out=80, in=180] (0.5, 1);
\draw [thick] (0.5,1) to [out=0, in=100] (0.9, 0.4);
\draw [thick] (0,0) to [out=-80, in=180] (0.5, -1);
\draw [thick] (0.5, -1) to [out=0, in=-100] (0.9, 0);
\node at (0,2.2) {$u^{-1}$};
\end{tikzpicture} \Bigg ) & = 2\;\;\begin{tikzpicture}[baseline={(current bounding box.center)},scale=0.3]
\draw [thick] (0,2) circle [radius=2];
\draw [thick] (-2,2) to (-0.6,2);
\draw [thick] (-0.8,2.2) to (-0.4,1.8);
\draw [thick] (-0.8, 1.8) to (-0.4, 2.2);
\draw [thick] (0,0) to [out=80, in=180] (0.5, 1);
\draw [thick] (0.5,1) to [out=0, in=100] (0.9, 0.4);
\draw [thick] (0,0) to [out=-80, in=180] (0.5, -1);
\draw [thick] (0.5, -1) to [out=0, in=-100] (0.9, 0);
\node at (0.7,2.2) {$u^{-1}$};
\end{tikzpicture} + \begin{tikzpicture}[baseline={(current bounding box.center)},scale=0.3]
\draw [thick] (0,2) circle [radius=2];
\draw [thick] (-1,0.2679) to (-0.6,2);
\draw [thick] (-0.8,2.2) to (-0.4,1.8);
\draw [thick] (-0.8, 1.8) to (-0.4, 2.2);
\draw [thick] (0,0) to [out=80, in=180] (0.5, 1);
\draw [thick] (0.5,1) to [out=0, in=100] (0.9, 0.4);
\draw [thick] (0,-1) to [out=180, in=-80] (-1, 0.2679);
\draw [thick] (0, -1) to [out=0, in=-100] (0.9, 0);
\node at (0.7,2.2) {$u^{-1}$};
\end{tikzpicture}\\
uD \Bigg ( \begin{tikzpicture}[baseline={(current bounding box.center)},scale=0.3]
\draw [thick] (0,2) circle [radius=2];
\draw [thick] (-2,2) to (-0.6,2);
\draw [thick] (-0.8,2.2) to (-0.4,1.8);
\draw [thick] (-0.8, 1.8) to (-0.4, 2.2);
\draw [thick] (-1.4142, 0.5858) to [out=40, in=140] (1, 0.5);
\draw [thick] (1.25, 0.3) to [out=-45, in=225] (2, 0.2);
\draw [thick] (2,0.2) to [out=45, in=-50] (1.732, 1);
\node at (0.7, 2.2) {$u^{-2}$};
\end{tikzpicture} \Bigg ) &= 6 \;\; \begin{tikzpicture}[baseline={(current bounding box.center)},scale=0.3]
\draw [thick] (0,2) circle [radius=2];
\draw [thick] (-1,0.2679) to (-0.6,2);
\draw [thick] (-0.8,2.2) to (-0.4,1.8);
\draw [thick] (-0.8, 1.8) to (-0.4, 2.2);
\draw [thick] (0,0) to [out=80, in=180] (0.5, 1);
\draw [thick] (0.5,1) to [out=0, in=100] (0.9, 0.4);
\draw [thick] (0,-1) to [out=180, in=-80] (-1, 0.2679);
\draw [thick] (0, -1) to [out=0, in=-100] (0.9, 0);
\node at (0.7,2.2) {$u^{-1}$};
\end{tikzpicture},
\end{align*}
and all other terms are zero.  Putting these together with the quantum
master equation we conclude that the coefficients $a$, $b$ satisfy  
\[\left\{ 
\begin{array}{c}
2b + \frac{1}{2} = 0\\ 
6a+b=0
\end{array}
\right. 
\]
Solving this gives the combinatorial model of the string vertex $S_{1,1}$:
\[ \iota(S_{1,1})=\frac{1}{24}\;\;\; \begin{tikzpicture}[baseline={(current bounding box.center)},scale=0.3]
\draw [thick] (0,2) circle [radius=2];
\draw [thick] (-2,2) to (-0.6,2);
\draw [thick] (-0.8,2.2) to (-0.4,1.8);
\draw [thick] (-0.8, 1.8) to (-0.4, 2.2);
\draw [thick] (-1.4142, 0.5858) to [out=40, in=140] (1, 0.5);
\draw [thick] (1.25, 0.3) to [out=-45, in=225] (2, 0.2);
\draw [thick] (2,0.2) to [out=45, in=-50] (1.732, 1);
\node at (0.7, 2.2) {$u^{-2}$};
\end{tikzpicture}
\;\;-\;\;\frac{1}{4}\;\;\; 
\begin{tikzpicture}[baseline={(current bounding box.center)},scale=0.3]
\draw [thick] (0,2) circle [radius=2];
\draw [thick] (0,0) to (0,1.4);
\draw [thick] (-0.2, 1.2) to (0.2, 1.6);
\draw [thick] (-0.2, 1.6) to (0.2, 1.2);
\draw [thick] (0,0) to [out=80, in=180] (0.5, 1);
\draw [thick] (0.5,1) to [out=0, in=100] (0.9, 0.4);
\draw [thick] (0,0) to [out=-80, in=180] (0.5, -1);
\draw [thick] (0.5, -1) to [out=0, in=-100] (0.9, 0);
\node at (0,2.2) {$u^{-1}$};
\end{tikzpicture}.\] 

\paragraph
This computation is essentially due to Costello (unpublished).  He
used it to compute the categorical Gromov-Witten invariant at $g=1$,
$n=1$ of a point (we take the $\cA_\infty$ algebra $A$ to be $\bbC$).
At $g=1$, $n=1$ the insertion of the $\psi$-class gives the coefficient
$1/24$ in the first term of $\iota(S_{1,1})$.  This agrees with the
geometric computation that
\[\int_{\overline{M}_{1,1}} \psi=\frac{1}{24}.\]

%% file: computer.tex
\section{The computer calculation}
\label{sec:computer}

In Section~\ref{sec:roadmap} we have argued that the computation of
Theorem~\ref{thm:mainthm} can be reduced to finding the correct chain
$\alpha$ which gives the lift $\tilde{\xi}$, and computing the 
values $T_2(\tilde{\xi}) = T_3(\tilde{\xi}) = 0$, 
\[ T_1(\tilde{\xi}) = \frac{1}{2}(\xi u^{-1})(\xi u^{-1}), \]
and 
\[ \langle \alpha, \xi \rangle_\Mukai = \frac{1}{12} E_2(\tau). \]
In this section we describe our initial approach to computing these
values, by  reducing the question of computing them to a large linear
algebra problem which can be solved by computer.  In the next section
we will give a purely mathematical deduction of Theorem~\ref{thm:mainthm}.

\paragraph
The fact that $T_3(\tilde{\xi}) = 0$ and 
\[ T_1(\tilde{\xi}) = \frac{1}{2}(\xi u^{-1})(\xi u^{-1}) \]
follows easily from the Kontsevich-Soibelman definition of the action
of ribbon graphs on cyclic chains of $A_\tau$.  For example the fact
that 
\[ \rho\big( 
\begin{tikzpicture}[baseline={(current bounding box.center)},scale=0.3]
\draw [thick] (0,2) circle [radius=2];
\draw [thick] (0,0) to (0,1.4);
\draw [thick] (-0.2, 1.2) to (0.2, 1.6);
\draw [thick] (-0.2, 1.6) to (0.2, 1.2);
\draw [thick] (0,0) to [out=80, in=180] (0.5, 1);
\draw [thick] (0.5,1) to [out=0, in=100] (0.9, 0.4);
\draw [thick] (0,0) to [out=-80, in=180] (0.5, -1);
\draw [thick] (0.5, -1) to [out=0, in=-100] (0.9, 0);
\end{tikzpicture}\big) (\xi) = 0 \]
follows form the fact that for any choice of basis vectors $x,y$ of
the algebra $A_\tau$ we have
\[ c_5(\xi, x, y, x^\chk, y^\chk) = 0, \]
where $x^\chk$, $y^\chk$ are dual basis vectors of $x, y$ with respect
to the pairing of the algebra $A_\tau$.  (This corresponds to the fact
that we need to label the leg of the graph labeled with a cross with
$\xi$, label the other half-edges of the graph by basis elements of
the algebra, and evaluate the corresponding cyclic operations at the
vertices and the duals of the pairing at the edges.)

Similarly, for the ribbon graph that appears in $S_{0,3}$ the input
edge is labeled $\xi$; the other two half edges adjacent to the
trivalent vertex will be labeled $\id_\cO$ in order to get a
non-trivial cyclic $c_3$, and then the output is read at the output
vertices as two copies of $\xi u^{-1}$.  
 
\paragraph
We next need to understand the conditions that single out the correct
lift 
\[ \tilde{\xi} = \xi + \alpha u + O(u^2) \in \HC_{-1}^-(A_\tau) \] 
of the class $\xi \in \HH_{-1}(A_\tau)$.  As mentioned before, for
degree reasons the only choice to be made is that of the chain
$\alpha \in CC_1(A_\tau)$.  This chain must satisfy
\[ b(\alpha) = -B(\xi) = -1\otimes\xi, \] 
and any choice of $\alpha$ satisfying this equation can be extended to
a class $[\tilde{\xi}]\in \HC_{-1}^-(A_\tau)$ which is uniquely
determined by $\alpha$.  Therefore we will talk about the lift
$[\tilde{\xi}]$ determined by $\alpha$. 

The following proposition will specify the correct choice of $\alpha$ to
match the geometric splitting in Lemma~\ref{lem:geomsplitting}.

\begin{Proposition}
\label{prop:charalpha}
Fix $\tau_0\in \bbH$ and consider the set of chains $\alpha\in CC_1(A_{\tau_0})$ satisfying
\[ b(\alpha) = -B(\xi) = -1\otimes\xi. \]
\begin{itemize}
\item[(i)] Among these chains there exists a unique one (up to
  $b$-exact chains) $\alpha^\GM$ with the property that the chain
  \[ b^{1|1}(\del_\tau\mu^*~|~ \alpha^\GM) \in CC_{-1}(A_{\tau_0}) \] 
  is $b$-exact, where $b^{1|1}$ is the operator defined
  in~\cite[3.11]{She} and $\del_\tau \mu^*$ is the derivative of the
  operations $\{\mu_k\}$, viewed as a Hochschild 2-cocycle in
  $\HH^2(A_{\tau_0})$. 

\item[(ii)] Let $[\tilde{\xi}]$ be the lift of $[\xi]$ to periodic
  cyclic homology corresponding to the unique chain in (i).  Then its
  image under the HKR isomorphism matches the class
  $[\tilde{\xi}]^\geom$ of Lemma~\ref{lem:geomsplitting}.
\end{itemize}
\end{Proposition}
\medskip

\begin{proof}
(i). In the proof of Lemma~\ref{lem:geomsplitting} we constructed a family
$[\xi]$ of classes in $\HH_{-1}(E_\tau)$ and a family of lifts
$[\tilde{\xi}]^\geom \in H^1_\dR(E_\tau)$, both parametrized by
$\tau\in \bbH$.  This family of lifts is uniquely characterized by the
local condition that it is Gauss-Manin flat.

Locally around $\tau_0\in\bbH$ pick representatives $\tilde{\xi}$ of
the images $[\tilde{\xi}]$ in $\HC_{-1}^-(A_{\tau_0})$ of the
geometric classes $[\tilde{\xi}]^\geom\in H^1_\dR(E_\tau)$ under the
inverse HKR isomorphism.  Each such $\tilde{\xi}$ will be of the form
\[ \tilde{\xi} = \xi + \alpha^\GM \cdot u + O(u^2) \] 
for some chain $\alpha^\GM\in CC_1(A_\tau)$ which depends on $\tau$.
Since $[\tilde{\xi}]$ depends smoothly on $\tau$ we can choose
$\alpha^\GM$ to also depend smoothly on $\tau$.

We have assumed that the geometric Gauss-Manin connection $\nabla^\GM$
is identified with the algebraic Getzler-Gauss-Manin connection 
$\nabla^\GGM$ under HKR.  The chain level Getzler-Gauss-Manin
connection applied to a family of chains $x\in CC_\sbt(A_\tau)$ is given by the
following formula, see~\cite[Definition 3.32]{She}:
\[ \nabla^\GGM_{\del_\tau}(x) = \del_\tau(x) -
  u^{-1} b^{1|1}(\del_\tau\mu^*~|~ x) -
  B^{1|1}(\del_\tau\mu^*~|~x). \]
Thus $\nabla^\GGM(\tilde{\xi})$ has the form
\[ \nabla^\GGM(\tilde{\xi}) = -b^{1|1}(\del_\tau\mu^*~|~\alpha^\GM) +
  O(u). \]
The flatness of the family $[\tilde{\xi}]^\geom$ implies that
$b^{1|1}(\del_\tau\mu^*~|~\alpha^\GM)$ is $b$-exact.

This shows that a family $\alpha$ with the required property exists,
at least locally around each $\tau_0$.  We will now show uniqueness.
Let $\alpha$ and $\alpha'$ be two chains satisfying the desired
properties.  Their difference is a Hochschild cycle:
\[ b(\alpha-\alpha') = -1\otimes \xi + 1\otimes \xi = 0, \]
so it gives a class $[\alpha-\alpha'] \in HH_1(A_\tau)$ for each
$\tau$.  The chain $\del_\tau \mu^* \in CC^2(A_{\tau_0})$ is a
cocycle, and its class in Hochschild cohomology is the negative of the
Kodaira-Spencer class
\[ \KS(\del_\tau) = -[\del_\tau \mu^*] \in \HH^2(A_{\tau_0}). ]\] 
(Note that the comparison between the algebraic and geometric
Kodaira-Spencer classes has a factor of (-1) built in;
see~\cite[3.32]{She}.)  Moreover, at homology level the class
\[ [ b^{1|1}(\del_\tau \mu^* ~|~ \alpha-\alpha') ] \in
  \HH_{-1}(A_{\tau_0}) \]
agrees with the contraction $-\KS(\del_\tau) \contract
[\alpha-\alpha']$.  
By assumption the former is $b$-exact.  

Therefore we have
\[ \KS(\del_\tau) \contract [\alpha-\alpha'] = 0 \]
in $\HH_{-1}(A_{\tau_0})$.  But for elliptic curves, contraction with
$\KS(\del_\tau)$ is an isomorphism $\HH_1(E_\tau)
\stackrel{\sim}{\lra} \HH_{-1}(E_\tau)$.  We conclude that
$[\alpha-\alpha'] = 0$ in $\HH_{-1}(E_\tau)$, and thus $\alpha^\GM$ is
unique up to $b$-exact chains.  (This uniqueness argument was suggested
to us by Sheridan.)

(ii).  Follows immediately from the arguments above for the proof of (i) and the
uniqueness of $[\tilde{\xi}]^\geom$.
\end{proof}

\paragraph
For a fixed value of $\tau\in\bbH$ the conditions 
\[ b(\alpha) = -1\otimes \xi \]
and 
\[ b^{1|1}(\del_\tau \mu^*~|~\alpha) = 0 \]
are linear in $\alpha$.  The structure constants of the
multiplications $\mu_k$ of $A_\tau$ can be computed using the formulas
in~\cite[1.2]{Pol}, and their derivatives with respect to $\tau$ can
be computed explicitly as well.  For example we have
\begin{align*}
g_{21} & = \frac{5}{6} g_{30} - g_{10}^2\\
g_{41} & = \frac{7}{10} g_{50} - 4g_{10}g_{30} \\
\del_\tau g_{10} & = \frac{g_{21}}{\twopii}-\frac{2g_{10}}{\tau-\taubar}.
\end{align*}
These are then used to write the linear operators $b$ and $b^{1|1}$ in
the basis for $A_\tau$ described in~(\ref{subsec:basis}).  The
resulting matrices allow us to express the equations
$b(\alpha) = -1\otimes \xi$, $b^{1|1}(\del_\tau\mu^*~|~\alpha) = 0$
and a solution $\alpha$ can be found in
$A_\tau^{\otimes 9} \oplus A_\tau^{\otimes 7} \oplus \cdots
A_\tau^{\otimes 3}$.  We then use $L^1$ optimization techniques to
find such a solution with a small number of non-zero terms; the best
such solution we found has 92 non-zero terms.

\paragraph
To complete the calculation of the potential we need to follow the
procedure outlined in~\cite[Section~11.6]{KonSoi} to compute the
result of inserting the various chains into the corresponding ribbon
graphs.

Evaluating the insertion of $\alpha$ is a large scale computation
which needs to be done by computer, but follows the same outline as
the ones above.   The results are 
\begin{align*}
& \rho \big(\begin{tikzpicture}[baseline={(current
                  bounding box.center)},scale=0.3] 
\draw [thick] (0,2) circle [radius=2];
\draw [thick] (-2,2) to (-0.6,2);
\draw [thick] (-0.8,2.2) to (-0.4,1.8);
\draw [thick] (-0.8, 1.8) to (-0.4, 2.2);
\draw [thick] (-1.4142, 0.5858) to [out=40, in=140] (1, 0.5);
\draw [thick] (1.25, 0.3) to [out=-45, in=225] (2, 0.2);
\draw [thick] (2,0.2) to [out=45, in=-50] (1.732, 1);
\end{tikzpicture}\big) (\alpha) = 0 \\
& \rho \big ( \begin{tikzpicture}[baseline={(current bounding box.center)},scale=0.3]
\draw [thick] (0,2) circle [radius=2];
\draw [thick] (-2,2) to (-0.6,2);
\draw [thick] (-0.8,2.2) to (-0.4,1.8);
\draw [thick] (-0.8, 1.8) to (-0.4, 2.2);
\draw [thick] (2,2) to (3.4,2);
\draw [thick] (3.2,2.2) to (3.6,1.8);
\draw [thick] (3.2, 1.8) to (3.6, 2.2);
\end{tikzpicture} \big ) (\alpha, \xi) = \frac{1}{12} E_2(\tau), 
\end{align*}
as required, completing the proof of Theorem~\ref{thm:mainthm}. 

\paragraph
\label{subsec:wrongsplit}
In the next section we will need to prove that the potential function
$F_{1,1}^{\mathsf B}(\tau)$ admits a finite limit at the cusp (i.e., as
$\tau\ra \mi\cdot\infty$).  Since the structure constants of
Polishchuk's algebra $A_\tau$ do extend to the cusp (in our basis --
see~\cite[2.5, Remark 2]{Pol}), the result would follow immediately if
we could argue that a choice for the chain $\alpha$ can be found which
also has a finite limit at the cusp.

Unfortunately we will not be able to prove directly the existence of such
a chain.  Instead we will use an auxiliary lift $[\tilde{\xi}']$ of
the family $[\xi]$ of Hochschild homology classes.  This lift will
{\em not} be Gauss-Manin flat, but by its very construction it will
have a finite limit at the cusp.  Using the lift $[\tilde{\xi}']$
instead of the correct, Gauss-Manin flat lift in the computations will
not compute the correct Gromov-Witten potential.  However, comparing
the result of this ``wrong'' computation with the correct one will
allow us to conclude that $F_{1,1}^{\mathsf B}(\tau)$ extends to the
cusp.

\paragraph
\label{subsec:alphaprime}
The desired family of chains
\[ \alpha' \in A_\tau^{\otimes 3} \oplus A_\tau^{\otimes 5} \oplus
  A_\tau^{\otimes 7} \]
which determines the lifts $\tilde{\xi}'$ can be written explicitly.  Taking
\begin{align*}
\alpha' = & \id_L  \theta  \eta  +\frac{1}{4} \cdot
            \id_\cO \theta  \eta +\frac{1}{2E_4(\tau)} \Big ( 9 E_2(\tau)  \cdot  \id_L\eta\xi\xi\theta\, +\\
& 60\cdot ( \theta\eta\xi\xi\theta\eta\xi +\eta\xi\theta\eta\xi\xi\theta )-12\cdot (\xi\theta\eta\xi\xi\theta\eta+\xi\xi\theta\eta\xi\theta\eta+\id_\cO\xi\theta\eta\xi\xi\xi) \\
& +36\cdot\id_\cO\xi\xi\theta\eta\xi\xi
  -24\cdot\id_\cO\xi\xi\xi\theta\eta\xi
  -60\cdot\id_L\eta\xi\xi\xi\xi\theta \Big ) 
\end{align*}
one can check by direct calculation that we have
\[ b(\alpha') = -1\otimes\xi = -B(\xi). \]
It is obvious from the above formulas that the chain $\alpha'$ extends
to the cusp: its coefficients have finite limit at $\tau = \mi\cdot\infty$.

%% file: computation.tex
\section{Proof of the main theorem}
\label{sec:computation}

In this final section we give a purely mathematical proof of
Theorem~\ref{thm:mainthm}, without relying on computer calculations.
The main idea is to carry out the same computation in two different
gauges, Polishchuk's modular gauge and the holomorphic gauge described
in Section~\ref{sec:gauge}, and with three different splittings of the
Hodge filtration.  Comparing the results of these computations will
allow us to determine the Gromov-Witten potential
$F_{1,1}^{\mathsf B}(\tau)$.

\paragraph
In a nutshell the computation of the categorical Gromov-Witten
invariant $F_{1,1}^{\mathsf B}(\tau)$ and of some of its variants can
be reduced to the following three steps:
\begin{itemize}
\item[(i)] Fix an $\cA_\infty$ model $A_\tau$ or $A_\tau^\hol$ of the
  derived category of the elliptic curve $E_\tau$. 
\item[(ii)] Construct a splitting of the Hodge
  filtration by finding a solution $\alpha$ of the equation
  \[ b(\alpha) = -B(\xi). \]
  (This $\alpha$ determines a lift $[\tilde{\xi}]$ of $[\xi]$ as in
  Section~\ref{sec:computer}.)
\item[(iii)] Compute the final invariant as
  \[ F_{1,1}^{\mathsf B}(\tau, \alpha) = T_2(\alpha) - \frac{1}{2} \langle
  \alpha, \xi\rangle_\Mukai,\]
where
\begin{align*}
T_2(\alpha) & = \frac{1}{24}\rho \big(\begin{tikzpicture}[baseline={(current
                  bounding box.center)},scale=0.3] 
\draw [thick] (0,2) circle [radius=2];
\draw [thick] (-2,2) to (-0.6,2);
\draw [thick] (-0.8,2.2) to (-0.4,1.8);
\draw [thick] (-0.8, 1.8) to (-0.4, 2.2);
\draw [thick] (-1.4142, 0.5858) to [out=40, in=140] (1, 0.5);
\draw [thick] (1.25, 0.3) to [out=-45, in=225] (2, 0.2);
\draw [thick] (2,0.2) to [out=45, in=-50] (1.732, 1);
\end{tikzpicture}\big) (\alpha)
\intertext{and}
\langle \alpha, \xi \rangle_\Mukai & = \rho \big ( 
\begin{tikzpicture}[baseline={(current bounding box.center)},scale=0.3]
\draw [thick] (0,2) circle [radius=2];
\draw [thick] (-2,2) to (-0.6,2);
\draw [thick] (-0.8,2.2) to (-0.4,1.8);
\draw [thick] (-0.8, 1.8) to (-0.4, 2.2);
\draw [thick] (2,2) to (3.4,2);
\draw [thick] (3.2,2.2) to (3.6,1.8);
\draw [thick] (3.2, 1.8) to (3.6, 2.2);
\end{tikzpicture} \big )(\alpha, \xi). 
\end{align*}
\end{itemize}
(Since we do not want to rely on computer calculations, we can not
assume that $T_2(\alpha) = 0$ as we have done before.)

In this section we will use three different choices for the chain
$\alpha$ in the above process, which give three different splittings
of the Hodge filtration and consequently three different potentials
(only one of which is of interest from the point of view of mirror
symmetry, see below).

\paragraph
The first choice of $\alpha$ is the chain $\alpha^\GM$ described in
Proposition~\ref{prop:charalpha}.  For this choice the corresponding
lifting $[\tilde{\xi}]$ is Gauss-Manin flat. The chain $\alpha^\GM$ is
characterized by the condition
\[ b^{1|1}(\del_\tau\mu^*~|~\alpha^\GM) = 0. \] 
Using it in the above procedure yields the correct Gromov-Witten
potential we are interested in,
\[F_{1,1}^{\mathsf B}(\tau) = F_{1,1}^{\mathsf B}(\tau,
  \alpha^\GM). \]

In a sense that will be made precise below the chain $\alpha^\GM$
depends in a holomorphic, but not modular fashion on $\tau$.  Thus
using it with the holomorphic $\cA_\infty$ model $A_\tau^\hol$ will
allow us to deduce that $F_{1,1}^{\mathsf B}(\tau)$ is holomorphic in
$\tau$.

\paragraph
The second choice of $\alpha$ is the 11-term chain $\alpha'$ described
in~(\ref{subsec:alphaprime}).  The resulting potential
$F_{1,1}^{\mathsf B}(\tau, \alpha')$ has the property that it extends
to the cusp.  Comparing it with $F_{1,1}^{\mathsf B}(\tau)$ will allow
us to deduce that this potential also admits a finite limit at the cusp.

\paragraph
\label{subsec:alphamod}
Finally, the third choice of chain $\alpha$ is a modular version
$\alpha^\modular$ of $\alpha^\GM$.  It is a solution of the equation
\[ b^{1|1}(\widehat{\del}_\tau\mu^*~|~\alpha^\GM) = 0, \] 
where $\widehat{\del}_\tau$ is the natural differential operator on
almost-holomorphic, modular forms described in~(\ref{subsec:delhat}).  Using it in
the modular model $A_\tau$ yields a new potential 
\[ F_{1,1}^{\mathsf B, \modular}(\tau) =  F_{1,1}^{\mathsf B}(\tau,
  \alpha^\modular). \]
By its very construction it will be obvious that $F_{1,1}^{\mathsf B,
  \modular}$ is modular.   

In Proposition~\ref{prop:calcdiff} we will argue that we have
\[ F_{1,1}^{\mathsf B}(\tau) - F_{1,1}^{\mathsf B, \modular}(\tau) =
  \frac{1}{4\pi\mi(\tau- \taubar)}. \]

The following theorem is a refinement of the statement of our main
result, Theorem~\ref{thm:mainthm}.

\begin{Theorem}
\label{thm:mainthm1}
The Gromov-Witten potential function $F_{1,1}^{\mathsf B}:\bbH\ra \bbC$ satisfies the
following properties:
\begin{enumerate}
\item[(a)] $F_{1,1}^{\mathsf B}$ is holomorphic.
\item[(b)] $F_{1,1}^{\mathsf B}$ extends to the cusp.
\item[(c)] The function
\[ F_{1,1}^{\mathsf B, \modular}(\tau)  = F_{1,1}^{\mathsf B}(\tau) - \frac{1}{4\pi\mi(\tau-\taubar)} \]
is modular of weight two. 
\end{enumerate}
Therefore we have
\[ F_{1,1}^{\mathsf B}(\tau) = -\frac{1}{24} E_2(\tau). \]
\end{Theorem}

\begin{Proof}
The proof of the fact that the potential $F_{1,1}^{\mathsf B}$
satisfies conditions (a)-(c) above is the content of
Propositions~\ref{prop:holomorphic},~\ref{prop:extends},
and~\ref{prop:calcdiff}.  The final statement of the theorem follows
from the fact that $-\frac{1}{24} E_2(\tau)$ is the unique function on
$\bbH$ that satisfies conditions (a)-(c).
\end{Proof}

\paragraph
In the remainder of this section we give precise arguments for the
fact that the potential $F_{1,1}^{\mathsf B}(\tau)$ satisfies
conditions (a)-(c) of Theorem~\ref{thm:mainthm1}.  We begin with a
lemma which shows that the term $T_2(\alpha)$ used in the definition
of the potential is not affected by changing $\alpha$ by any multiple
of the class $[\Omega] \in HH_1(A_\tau)$ defined in
Proposition~\ref{prop:defOmegaxi}.

\begin{Lemma}
\label{lem:insertomega}
We have
\[ \rho \Big( \begin{tikzpicture}[baseline={(current bounding box.center)},scale=0.3]
\draw [thick] (0,2) circle [radius=2];
\draw [thick] (-2,2) to (-0.6,2);
\draw [thick] (-0.8,2.2) to (-0.4,1.8);
\draw [thick] (-0.8, 1.8) to (-0.4, 2.2);
\draw [thick] (-1.4142, 0.5858) to [out=40, in=140] (1, 0.5);
\draw [thick] (1.25, 0.3) to [out=-45, in=225] (2, 0.2);
\draw [thick] (2,0.2) to [out=45, in=-50] (1.732, 1);
\end{tikzpicture} \Big) ([\Omega]) =0.\]
\end{Lemma}

\begin{Proof}
The Calabi-Yau structure of $A_\tau$ gives an identification 
\[ HH^\sbt(A_\tau) \cong HH_{1-\sbt}(A_\tau). \] 
This turns the former into a graded commutative Frobenius
algebra with pairing given by the Mukai pairing on homology.  In turn
this yields a $2$-dimensional topological field theory.  

Moreover, under this identification the class $[\Omega] \in HH_1(A_\tau)$
is matched with a multiple of $1\in HH^0(A_\tau)$.  The field theory
interpretation implies that
\[\rho\Big(\begin{tikzpicture}[baseline={(current bounding box.center)},scale=0.3]
\draw [thick] (0,2) circle [radius=2];
\draw [thick] (-2,2) to (-0.6,2);
\draw [thick] (-0.8,2.2) to (-0.4,1.8);
\draw [thick] (-0.8, 1.8) to (-0.4, 2.2);
\draw [thick] (-1.4142, 0.5858) to [out=40, in=140] (1, 0.5);
\draw [thick] (1.25, 0.3) to [out=-45, in=225] (2, 0.2);
\draw [thick] (2,0.2) to [out=45, in=-50] (1.732, 1);
\end{tikzpicture}\Big) ([\Omega])= \mbox{ a multiple of }\Tr_{HH^\sbt(A_\tau)}(1). \]
The term $\Tr_{HH^\sbt(A_\tau)}(1)$ is the supertrace of multiplication
by $1$ in the Hochschild cohomology ring $\HH^\sbt(A_\tau)$, which is
zero (the Euler characteristic of the graded vector space
$\HH^\sbt(A_\tau)$). 
\end{Proof}

\paragraph
Let $\alpha_1(\tau), \alpha_2(\tau)$ be two families of solutions of
the equation $b(\alpha) = -B(\xi)$.  Then for every $\tau$ the chain
$\alpha_1(\tau)-\alpha_2(\tau)$ is $b$-closed, and therefore it
defines a class $[\alpha_1(\tau)-\alpha_2(\tau)] \in \HH_1(A_\tau)$ which must be
a multiple $c(\tau)\cdot [\Omega]$ of $[\Omega]$.

\begin{Corollary}
\label{cor:diffalpha}
The difference between the potentials obtained from $\alpha_1$ and
$\alpha_2$ is 
\[ F_{1,1}^{\mathsf B}(\tau, \alpha_1) - F_{1,1}^{\mathsf B}(\tau,
  \alpha_2) = -\frac{c(\tau)}{2}. \]
\end{Corollary}

\begin{Proof}
We have
\begin{align*}
F_{1,1}^{\mathsf B}(\tau, \alpha_1) - F_{1,1}^{\mathsf B}(\tau,
  \alpha_2) & = T_2(\alpha_1-\alpha_2) - \frac{1}{2}\langle
              \alpha_1-\alpha_2, \xi\rangle_\Mukai \\ 
& = T_2(c(\tau)\cdot [\Omega]) - \frac{1}{2} \langle c(\tau)\cdot
  [\Omega], [\xi]\rangle_\Mukai \\
& = -\frac{c(\tau)}{2}. 
\end{align*}
\end{Proof}

\begin{Proposition}
\label{prop:holomorphic}
The categorical Gromov-Witten invariant $F_{1,1}^{\mathsf B}(\tau)$ is
a holomorphic function on $\bbH$.  Moreover, it is an element of
weight two in the graded ring $\tM(\Gamma)_{(0)}$ obtained by
localizing $\tM(\Gamma)$ at the multiplicative set of elements of
homogeneous weight.
\end{Proposition}

\begin{Proof}
By~(\ref{subsec:invariance}) and Theorem~\ref{thm:holgauge} it will be
enough to prove the statement when using the algebra $A_\tau^\hol$
that we defined in~(\ref{subsec:defatauhol}) as the chosen
$\cA_\infty$ model of the derived category of the elliptic curve
$E_\tau$.  Therefore we fix the algebra $A = A_\tau^\hol$.
 
In Proposition~\ref{prop:charalpha} we have argued that the equations
defining the chain $\alpha^\GM$ can always be solved locally around
every $\tau_0\in \bbH$, and the solution is unique up to $b$-exact
chains.  Moreover, the corresponding lift of $[\xi]$ will be
Gauss-Manin-flat.  The proof of this result did not depend on the
specific algebra we used, only on the fact that the section $\xi$ of
the family $\{A_\tau\}_{\tau\in \bbH}$ was flat with respect to
differentiation by $\del_\tau$.  Using any local choice
$\alpha^{\GM}$ of such $\alpha$ allows us to define
$F_{1,1}^{\mathsf B}$ locally on $\bbH$.

The action of ribbon graphs on $CC_\sbt(A)$ is a chain map, in the
sense that it intertwines the $\del$ operator on ribbon graphs and the
$b$ operator on chains.  Both ribbon graphs involved in the definition
of the potential are trivalent, therefore they are sent to zero by
$\del$.  This implies that modifying $\alpha^{\GM}$ by a $b$-exact chain
does not change the result of its insertion in these two graphs.  We
conclude that the locally defined potentials glue to a well-defined
global potential on the entire upper half plane $\bbH$. 

The structure constants of the algebra $A_\tau^\hol$ are holomorphic,
quasi-modular forms in the ring $\tM(\Gamma)$.  All the operations
used to compute the Gromov-Witten potential (differentiating, solving
linear systems, inserting chains into ribbon graphs) always keep us
inside $\tM(\Gamma)_{(0)}$.  (We need to allow for divisions by
homogeneous weight elements of $\tM(\Gamma)$ in order to solve linear
systems.)  We conclude that $F_{1,1}^{\mathsf B}(\tau)$ is an element
of $\tM(\Gamma)_{(0)}$, and an easy computation shows that its weight
is two. 

We have argued above that $F_{1,1}^{\mathsf B}$ is defined everywhere
on $\bbH$.  Checking that it is holomorphic is a local computation, which
can be done with one local choice of $\alpha^{\GM}$.
\end{Proof}

\begin{Proposition}
\label{prop:extends}
The potential $F_{1,1}^{\mathsf B}(\tau)$ extends to the cusp: the
limit
\[ \lim_{\tau\ra \mi\cdot\infty} F_{1,1}^{\mathsf B}(\tau) \]
exists and is finite.
\end{Proposition}

\begin{Proof} 
We use the modular gauge $A_\tau$.  The structure constants of this
family of $\cA_\infty$ algebras extend to the cusp~\cite[2.5, Remark
2]{Pol}.  The chain $\alpha'$ defined in~(\ref{subsec:alphaprime})
also extends to the cusp.  The terms $T_2(\alpha')$ and
$\langle \alpha', \xi\rangle_\Mukai$ are (very complicated) polynomial
expressions in the structure constants of the algebra and the
coefficients of $\alpha'$.  Therefore the potential
$F_{1,1}^{\mathsf B}(\tau, \alpha')$ obtained from the chain $\alpha'$
extends to the cusp.

The chains $\alpha'$ and $\alpha^\GM$ both satisfy
$b(\alpha) = -B(\xi)$, thus we are in the setup of
Corollary~\ref{cor:diffalpha}.  If we argue that the function
$c(\tau)$ defined by the equality
\[ [\alpha^\GM-\alpha'] = c(\tau)\cdot[\Omega] \]
extends to the cusp, we will be able to conclude that 
\[ F_{1,1}^{\mathsf B}(\tau) = F_{1,1}^{\mathsf B}(\tau, \alpha') -
  \frac{c(\tau)}{2} \]
also extends to the cusp, which is what we need to prove.

The expression $b^{1|1}(\del_\tau \mu^*~|~\alpha')$ is a polynomial
expression in the structure constants of $A_\tau$, their derivatives,
and the coefficients of $\alpha'$, all of which have limits at the
cusp; therefore this expression also extends to the cusp.  Using the
defining property of $\alpha^\GM$ that
\[ b^{1|1}(\del_\tau \mu^*~|~\alpha^\GM) = 0 \]
we conclude that $b^{1|1}(\del_\tau \mu^*~|~\alpha'-\alpha^\GM)$ also
extends to the cusp. 

By HKR we conclude that 
\[ -\KS(\del_\tau) \contract c(\tau)\cdot[\Omega] \]
extends to the cusp.  The class $\KS(\del_\tau)$ was computed
in~(\ref{subsec:KS}), 
\[ \KS(\del_\tau) = -\frac{1}{\tau-\taubar} \frac{\del}{\del
    z}\dzbar \] 
and therefore
\[ -\KS(\del_\tau) \contract [\Omega] = \frac{1}{\tau-\taubar}
  \frac{\del}{\del z} \dzbar \contract [\twopii\cdot dz] =
  \frac{\twopii}{\tau-\taubar}\dzbar = \twopii \cdot \xi. \]
It follows that $c(\tau)$ admits a finite limit at the cusp. 
\end{Proof}

\paragraph
From now on we will position ourselves in the
modular gauge -- we take Polishchuk's algebra $A_\tau$ as the
$\cA_\infty$ model for the elliptic curve.

Fix any $\tau\in \bbH$.  We have already shown that the equation
\[ b^{1|1}(\del_\tau \mu^*~|~\alpha) = 0 \] 
along with $b(\alpha) = -B(\xi)$ defines a unique (up to $b$-exact
terms) chain $\alpha^\GM$ which we used to define the potential
$F_{1,1}^{\mathsf B}(\tau)$.

However, the resulting chain $\alpha^\GM$ is not modular.  The point where
modularity breaks down is in the computation of $\del_\tau \mu^*$: even
though the coefficients of $\mu^*$ are in the ring
$\widehat{M}(\Gamma)$ of almost holomorphic modular forms, the
coefficients of $\del_\tau \mu^*$ are no longer modular.

\paragraph
This observation suggests how to modify the above computation so that
all the intermediate computations (and the result) stay in
$\widehat{M}(\Gamma)$.  Namely, if we replace the equation above with
\[ b^{1|1}(\widehat{\del}_\tau \mu^*~|~\alpha^\modular) = 0 \] 
where $\widehat{\del}_\tau$ is the natural differential operator on
$\widehat{M}(\Gamma)$
\[ \widehat{\del}_\tau = \del_\tau + \frac{\mathsf{wt}}{\tau-\taubar} \] 
the solution $\alpha^\modular$ will depend in a modular, almost
holomorphic way on $\tau$.  

Note that locally around every point in $\bbH$ such a solution
$\alpha^\modular$ exists and is unique up to a $b$-exact chain because
the operator $\KZ$ of Theorem~\ref{thm:KZ} induces an isomorphism
between solutions $\alpha^\modular$ of the above equation and
solutions $\alpha$ of the original equation
\[ b^{1|1}(\del_\tau\mu^*~|~\alpha). \]

\paragraph
We use the solution $\alpha^\modular$ to define a new potential 
\[ F_{1,1}^{\mathsf B, \modular}(\tau) = F_{1,1}^{\mathsf B}(\tau, 
  \alpha^\modular). \]
As the operator $\KZ$ is a differential ring isomorphism, it follows
that we have
\[ F_{1,1}^{\mathsf B, \modular} = \KZ(F_{1,1}^{\mathsf B}), \]
and in particular $F_{1,1}^{\mathsf B, \modular}$ is a modular form in
$\widehat{M}(\Gamma)_{(0)}$ of weight two. 

\begin{Proposition}
\label{prop:calcdiff}
We have
\[ F_{1,1}^{\mathsf B}(\tau) - F_{1,1}^{\mathsf B, \modular} =
  \frac{1}{4\pi\mi(\tau-\taubar)}. \] 
Therefore $F_{1,1}^{\mathsf B}(\tau)$ satisfies condition (c) of
Theorem~\ref{thm:mainthm1}.  
\end{Proposition}

\begin{Proof}
Consider the difference 
\[ \delta = \alpha^\GM-\alpha^\modular. \] 
It is a $b$-closed chain in $CC_1(A_\tau)$, thus it gives a class
$[\delta]$ in $\HH_1(A_\tau)$.   We will show that 
\[ [\delta]  = -\frac{1}{\twopii(\tau-\taubar)} [\Omega]. \]
Corollary~\ref{cor:diffalpha} then implies the result. 

The chains $\alpha^\GM$ and $\alpha^\modular$ define lifts of $\xi$ which
we denote by $\tilde{\xi}$ and $\tilde{\xi}^\modular$, respectively.
The former is flat with respect to the Gauss-Manin connection.  The
main idea is to calculate
\[ \nabla^\GM_{\del_\tau}([\xi]^\modular) \]
in two different ways, and to compare the results.

The definition of the Getzler-Gauss-Manin connection 
\[ \nabla^\GGM_{\del_\tau}(x) = \del_\tau(x) - u^{-1}
  b^{1|1}(\del_\tau\mu^*~|~ x) - B^{1|1}(\del_\tau\mu^*~|~x) \]
implicitly assumes that we have chosen a connection on the family
$A_\tau$ of $\cA_\infty$ algebras over $\bbH$ with the property that
the basis elements of $A_\tau$ are flat for this connection.  This
allows us to write expressions like $\del_\tau(x)$, by which we mean
that we write $x$ in this chosen basis and differentiate the
coefficients. 

However there is no reason to insist on the use of such a connection
-- {\em any} connection $\nabla$ will work to define the
Getzler-Gauss-Manin connection, at least at the level of homology.  We
just need to replace the formula above by
\[ \nabla^\GGM_{\del_\tau}(x) = \nabla_{\del_\tau} (x) - u^{-1}
  b^{1|1}(\nabla_{\del_\tau}\mu^*~|~ x) -
  B^{1|1}(\nabla_{\del_\tau}\mu^*~|~x), \]
where now we think of applying $\nabla_{\del_\tau}$ to arbitrary tensors.

We will use two different connections $\nabla^{\mathsf{std}}$ and
$\nabla^\modular$ on the bundle $A_\tau$ over $\bbH$.  The first one,
$\nabla^{\mathsf{std}}$, is the standard connection for which the
basis vectors $\id_{\cO}, \id_L, \theta, \eta, \xi, \xi_L$ are flat.
Computing $\nabla^\GGM{\del_\tau}(\tilde{\xi})$ using $\nabla^{\mathsf{std}}$
yields
\begin{align*}
\nabla^\GGM_{\del_\tau}(\tilde{\xi}^\modular) & = \nabla^{\mathsf{std}}_{\del_\tau}(\xi) -
  b^{1|1}(\nabla^{\mathsf{std}}_{\del_\tau}\mu^*~|~\alpha^\modular) + O(u) \\
& = -b^{1|1}(\del_\tau\mu^*~|~\alpha^\modular) + O(u).
\end{align*}

Define a second connection $\nabla^\modular$ on the bundle $A_\tau$ of
$\cA_\infty$ algebras by the formula
\[ \nabla^\modular_{\del_\tau} (x) = -\frac{\wt(x)}{\tau-\taubar}x \] 
for any basis vector $x$ of $A_\tau$, where the weight of basis
elements was defined  in~(\ref{subsec:wtbasis}).  A straightforward
computation shows that 
\[ \nabla_{\del_\tau}^\modular \mu_k = \widehat{\del}_\tau \mu_k, \]
where the left hand side refers to the differentiation of the tensor
$\mu_k$ with respect to the connection $\nabla$, while the right hand
side refers to the differentiation of the structure constants of the
same tensor using the operator $\widehat{\del}_\tau$. 

Computing  $\nabla^\GGM_{\del_\tau}(\tilde{\xi}^\modular)$ using the connection
$\nabla^\modular$ yields 
\begin{align*}
\nabla^\GGM_{\del_\tau}(\tilde{\xi}^\modular) & =
\nabla_{\del_\tau}^\modular(\xi) - b^{1|1}(\nabla_{\del_\tau}^\modular\mu^*~|~\alpha^\modular) + O(u)\\
& = -\frac{1}{\tau-\taubar}\xi -
  b^{1|1}(\widehat{\del}_\tau\mu^*~|~\alpha^\modular) + O(u)\\
& = -\frac{1}{\tau-\taubar}\xi +O(u).
\end{align*}

Equating the two calculations above yields
\[ b^{1|1}(\del_\tau\mu^*~|~\delta)  =
  b^{1|1}(\del_\tau\mu^*~|~\alpha^\GM-\alpha^\modular)=-\frac{1}{\tau-\taubar}\xi + O(u). \]
Passing to Hochschild homology we get
\[ [\del_\tau \mu^*] \contract [\delta] = -\frac{1}{\tau-\taubar}[\xi]. \]
Since $[\delta] = c(\tau)\cdot[\Omega]$ for some $c(\tau)$, under the
HKR isomorphism this becomes 
\[ -\KS(\del_\tau) \contract (c(\tau)\cdot [\Omega]) = -\twopii\cdot c(\tau)\cdot
  \KS(\del_\tau) \contract [dz] =
  -\frac{1}{(\tau-\taubar)^2}\dzbar, \]
which forces $c(\tau)$ to equal
\[ c(\tau) = -\frac{1}{\twopii(\tau-\taubar)}. \]
\end{Proof}

\paragraph {\bf Remark.}
It is worth noting that the modular splitting obtained from the chain
$\alpha^\modular$ corresponds to the complex conjugate splitting from
complex algebraic geometry.  This follows from a direct computation.